\documentclass[opre,nonblindrev]{informs3}

\OneAndAHalfSpacedXI 
\usepackage{endnotes}
\let\footnote=\endnote

\usepackage{mathtools}
\usepackage{algorithm}
\usepackage{algpseudocode}
\usepackage{booktabs}

\usepackage[edges]{forest}

\DeclareMathOperator{\adj}{adj}

\DeclareMathOperator{\projone}{PROJ_1}
\DeclareMathOperator{\projtwo}{PROJ_2}
\DeclareMathOperator{\projh}{PROJH}
\DeclareMathOperator{\psc}{psc}
\DeclareMathOperator{\pscset}{PSC}
\DeclareMathOperator{\red}{red}
\DeclareMathOperator{\redset}{RED}
\DeclareMathOperator{\ldt}{ldt}
\DeclareMathOperator{\ldcf}{ldcf}
\DeclareMathOperator{\sylvha}{SylvHa}

\newcommand{\mat}[1]{\boldsymbol{\mathbf{#1}}}

\usepackage{graphicx}
\usepackage{multirow}
\usepackage{hhline}

\usepackage[small, margin=1cm]{caption}

\usepackage{appendix}
\usepackage{color}
\definecolor{strcolor}{rgb}{0.6, 0.2, 0.6}
\definecolor{commentcolor}{rgb}{0.3125, 0.5, 0.3125}
\definecolor{keycol}{rgb}{0, 0, 1}

\usepackage{amssymb}
\usepackage{amsmath}
\usepackage{bbm}

\usepackage{url}

\newcommand {\bea}{\begin{eqnarray}}
	\newcommand {\eea}{\end{eqnarray}}

\def\blot{\quad \mbox{$\vcenter{ \vbox{ \hrule height.4pt
				\hbox{\vrule width.4pt height.9ex \kern.9ex \vrule width.4pt}
				\hrule height.4pt}}$}}

\usepackage{natbib}
\bibpunct[, ]{(}{)}{,}{a}{}{,}%
\def\bibfont{\small}%
\TheoremsNumberedThrough     
\ECRepeatTheorems

\EquationsNumberedThrough    

\gdef\AQ#1{}
\gdef\CQ#1{}

\begin{document}
	
\def\COPYRIGHTHOLDER{INFORMS}%
\def\COPYRIGHTYEAR{2017}%
\def\DOI{\fontsize{7.5}{9.5}\selectfont\sf\bfseries\noindent https://doi.org/10.1287/opre.2017.1714\CQ{Word count = 9740}}

	\RUNAUTHOR{Chan et~al.} %

	\RUNTITLE{Exact sensitivity analysis of Markov reward processes via algebraic geometry}

\TITLE{Exact sensitivity analysis of Markov reward processes via algebraic geometry}


	\ARTICLEAUTHORS{

\AUTHOR{Timothy C. Y. Chan, Muhammad Maaz}
\AFF{Department of Mechanical and Industrial Engineering,
University of Toronto}


}
	

	\ABSTRACT{We introduce a new approach for deterministic sensitivity analysis of Markov reward processes, commonly used in cost-effectiveness analyses, via reformulation into a polynomial system. Our approach leverages cylindrical algebraic decomposition (CAD), a technique arising from algebraic geometry that provides an exact description of all solutions to a polynomial system. While it is typically intractable to build a CAD for systems with more than a few variables, we show that a special class of polynomial systems, which includes the polynomials arising from Markov reward processes, can be analyzed much more tractably. We establish several theoretical results about such systems and develop a specialized algorithm to construct their CAD, which allows us to perform exact, multi-way sensitivity analysis for common health economic analyses. We develop an open-source software package that implements our algorithm. Finally, we apply it to two case studies, one with synthetic data and one that re-analyzes a previous cost-effectiveness analysis from the literature, demonstrating advantages of our approach over standard techniques. Our software and code are available at: \url{https://github.com/mmaaz-git/markovag}.
    }





	
	%
	
\maketitle

\section{Introduction}

This paper develops a new approach to conduct exact, deterministic sensitivity analyses of a Markov reward process, motivated by their ubiquitous use in cost-effectiveness analyses (CEA). Our approach leverages ideas from algebraic geometry, particularly cylindrical algebraic decomposition (CAD), and applies them to the mathematical structure of common cost-effectiveness analyses, which enables exact analysis more efficiently than for general polynomial systems. The implication is that exact multi-way (where multiple parameters are varied simultaneously) sensitivity analysis is possible, which allows a policymaker to fully describe complex parameter regimes where new technologies or medical interventions are cost-effective.

A Markov reward process often forms the modeling foundation for a CEA. For example, in a health economic analysis, states may represent patient health states and the Markov process would model disease progression of the patient as lifetime costs and benefits accumulate \citep{rudmik2013health}. A new medical intervention would be evaluated based on how it alters the primitives of the Markov reward process. For example, a promising intervention might lower the transition probability to worse health states or increase rewards associated with quality of life in healthier states. After calculating lifetime benefits and costs using the model, a policymaker may be interested in whether the total reward associated with the intervention exceeds a given threshold, or whether the reward of one intervention exceeds that of an alternative. Other metrics of interest in a CEA may include the net monetary benefit (benefit multiplied by willingness-to-pay per unit benefit minus cost) or the incremental cost-effectiveness ratio (difference in costs of two interventions divided by their difference in benefits).

All of the above metrics depend on the parameters of the Markov reward process. However, these parameter values may be subject to significant uncertainty, especially when the intervention being evaluated is new and there is limited data about it. Often, values are drawn from related empirical studies in the literature or simply assumed based on expert opinion. Thus, sensitivity analysis is critical for any CEA. A typical approach is a one-way deterministic sensitivity analysis, which means that one parameter is varied within a range or set of values while all others are fixed at some nominal value. This approach is straightforward but does not capture the joint effect of multiple parameters. Multi-way sensitivity analyses create a multi-dimensional grid of test points over the parameter space, which measures the impact of parameter interactions, but quickly becomes intractable in the number of parameters tested and grid granularity. In a systematic review of CEAs, \citet{jain2011sensitivity} found that 86\% of studies in their sample conducted a one-way sensitivity analysis, but only 45\% conducted a multi-way sensitivity analysis, likely owing to these difficulties. In practice, multi-way sensitivity analyses rarely extend beyond two parameters \citep{briggs1994uncertainty}. We note that there is another type of sensitivity analysis known as probabilistic sensitivity analysis, where parameter values are drawn from distributions and their joint effect simulated \citep{baio2015probabilistic}. However, they face a similar issue as parameter distributions are often arbitrarily chosen. In contrast, our method makes no assumptions on the ranges nor distributions of the parameters, and instead enumerates the full range of parameters that yields the desired result.

We make the key observation that the questions typically asked in a cost-effectiveness analysis based on a Markov reward process can be described as a system of polynomial inequalities. Determining whether an intervention remains cost-effective if parameter values vary within given intervals is thus equivalent to determining whether a polynomial system satisfies a set of inequalities over those parameter intervals. Hence, we study deterministic sensitivity analysis through the lens of algebraic geometry, which provides tools that facilitate analysis of multivariate polynomial systems, such as \textit{cylindrical algebraic decomposition} \citep{collins1974quantifier}. CAD was the first practical algorithm for solving systems of polynomial inequalities and works by decomposing the multidimensional real space into cylindrical cells over which each polynomial is sign invariant. Once these cells are defined, the algorithm can easily check feasibility of the polynomial system over each cell \citep{basu}. In doing so, we obtain a tree-like representation of the whole space over which the polynomial system holds. In this paper, we will use CAD to fully represent the cost-effectiveness region in a multi-way sensitivity analysis. While analyzing general polynomial systems using CAD remains computationally challenging, the polynomial systems of interest in a Markov reward process-based CEA can be analyzed much more tractably.

Our main contributions are as follows.

\begin{enumerate}
    \item \textit{Semialgebraic representation of Markov reward processes.} We show that sensitivity analysis of common cost-effectiveness analysis quantities, including bounding or comparing total benefits, total costs, incremental cost-effectiveness ratios, and net monetary benefits, is equivalent to determining whether certain polynomial systems are feasible. Thus, describing the cost-effective parameter space can be done by analyzing the CAD of these systems.
    \item \textit{Cylindrical algebraic decomposition.} We demonstrate that the polynomial systems induced by the aforementioned analyses belong to a special class that makes the CAD construction more efficient. This class accommodates important considerations such as the transition probability matrix possessing the increasing failure rate property. We develop a specialized version of the general CAD algorithm for this class of systems. We show that this CAD has a singly exponential size, compared to the doubly exponential size in a general CAD. 
    \item \textit{Software.} We develop a Python package, \texttt{markovag}, that implements our algorithms. It can construct polynomial systems representing sensitivity analysis of common CEA metrics, analytically characterize the boundary in a multi-way sensitivity analysis, and construct the CAD. 
    \item \textit{Case studies.} To demonstrate the CAD approach to sensitivity analysis, we apply our algorithms and software in  two case studies. The first case study uses synthetic data to show that a traditional parameter grid search could easily mischaracterize a non-linear cost-effectiveness boundary. In the second case study, we re-analyze a real CEA from the literature and show that our approach reveals a larger cost-effective parameter space than in the original analysis and elucidates relationships between model parameters that would not be otherwise obvious.
\end{enumerate}


\section{Literature Review}

This paper relates to three different bodies of literature, across the domains of health economics, stochastic processes, and algebraic geometry.

\subsection{Cost-effectiveness analysis}

Markov models have been a ubiquitous tool for cost-effectiveness analysis in healthcare \citep{carta2020use}. They are part of the larger class of multistate models \citep{hougaard1999multi} that model transitions between health states. Markov models in healthcare are typically analyzed either using linear algebra or via cohort simulations \citep{sonnenberg1993markov}. As these models exhibit parameter uncertainty, sensitivity analysis is an essential step \citep{briggs1994uncertainty, jain2011sensitivity, rudmik2013health}. Indeed, it is recommended by health economics professional societies \citep{briggs2012model} and even mandated by policymakers \citep{andronis2009sensitivity}. 
However, when sensitivity analyses are limited to one or two parameters, larger interaction effects or correlations may be missed \citep{vreman2021application}. Furthermore, the range of parameter values to test are often chosen arbitrarily. In contrast, our method can perform arbitrary multi-way sensitivity analyses, uncovering interactions between parameters, and will identify the full range of parameter values that lead to a cost-effectiveness result holding.

\subsection{Uncertainty in Markov models}

The notion of Markov chains with imprecise parameters is well-studied \citep{caswell2019sensitivity, dai1996sensitivity, caswell2013sensitivity, hermans2012characterisation}. For example, \citet{de2014sensitivity} develop a generalization of the Perron-Frobenius Theorem when the transition probabilities lie within a credal set. \citet{blanc2008markov} study Markov chains with row-wise uncertainty to compute bounds on hitting times and stationary distributions. However, this line of research focuses on Markov chains, not Markov reward processes, and as a result is silent on derived metrics arising from cost-effectiveness analysis. There is also a large literature on Markov decision processes with uncertain parameters \citep{el2005robust, iyengar2005robust, delage2010percentile, wiesemann2013robust, goyal2023robust, grand2024convex}. 

The paper closest to ours is that of \citet{goh}, who develop a method for finding the maximum or minimum infinite horizon total reward for a Markov reward process, where the set of probability matrices has a row-wise structure. They show that such a problem can be formulated as a Markov decision process, and solving it using standard policy iteration provides the optimal solution to the original problem. The connection to CEA is that if a policymaker requires an outcome metric to be above a certain threshold, but the maximum value of that metric lies below the threshold, then they should reject the treatment. However, simply knowing the extrema would be insufficient if the threshold lies between the maximum and minimum, as the policymaker would not know over which parameter values the threshold is met. In contrast, our work allows a policymaker to know exactly which parameter values attain or violate a given inequality. 

\subsection{Cylindrical algebraic decomposition}

A cylindrical algebraic decomposition (CAD) is a tree-like decomposition of the real space representing the solutions to a polynomial system, which is how we apply it in this paper. The algorithm to construct a CAD tree was initially developed by \citet{collins1974quantifier} to solve the problem of the existential theory of the reals, which asks for a satisfying assignment of real numbers to a Boolean combination of polynomials. This problem was known to be solvable theoretically by \cite{tarski1951decision} but it was not until \citet{collins1974quantifier} that a practical algorithm was developed. This problem is important in computer science as it fits into the general framework of satisfiability modulo theories (SMT), a generalization of the classical Boolean satisfiability problem (SAT) \citep{papa} to statements with variables that can be numbers or even data structures \citep{de2011satisfiability}. Solving polynomial systems with CAD has long been used in robotics motion planning \citep{canny1987new, schwartz1990algorithmic}.

More broadly, there are some strong connections between techniques from algebraic geometry and the operations research literature, particularly for polynomial optimization and semidefinite programming \citep{lasserre2001global,parrilo2003semidefinite, blekherman2012semidefinite, parrilo2019sum}. However, to date, techniques from this field have not been applied to the study of Markov reward processes. Our technical contribution will be to apply concepts from algebraic geometry to study the solutions to our class of polynomial systems, derived from Markov reward process sensitivity analysis.


The computer algebra literature has identified special cases of polynomial systems for which CAD can be simplified. There is a recent body of work on speeding up CAD when the system contains equalities \citep{england2020cylindrical}, which do arise in the class of systems we study, but the algorithm remains asymptotically doubly exponential. Incremental CAD \citep{kremer2020fully} extends a CAD of a polynomial system to a CAD of the same system with an additional inequality. We use a similar idea of extending a CAD by incrementally adding inequalities in our algorithm. Lastly, \citet{strzebonski2010computation} developed an algorithm to construct a CAD from Boolean formulas of CADs. Our contribution to computer algebra is the identification of a special but broad class of polynomial systems under which we can construct the CAD more tractably.

\section{Semialgebraic Representations of Markov Reward Processes}

We consider a discrete-time discounted Markov reward process characterized by a Markov chain with $n$ states, a transition probability matrix $\mat{P} \in \mathbb{R}^{n \times n}$ ($p_{ij}$ is the transition probability of going from state $i$ to $j$), a reward vector $\mat{r} \in \mathbb{R}^n$ ($r_i$ is the reward for being in state $i$ for one period), an initial state distribution $\mat{\pi} \in \mathbb{R}^n$ ($\pi_i$ is the probability the process starts in state $i$), and a discount rate $\lambda \in (0,1)$. Let $\mat{I}$ denote the identity matrix.

The expected discounted reward over a finite horizon of length $t$, $R_t$, or over an infinite horizon, $R_\infty$, can be written as \citep{putermanmarkov}: 
%
\begin{equation}
    \label{eq:finite_reward}
    R_t := \sum_{m=0}^t \mat{\pi}^\top \lambda^m \mat{P}^m \mat{r}
\end{equation}
and
\begin{equation}
    \label{eq:infty_reward}
    R_\infty := \mat{\pi}^\top (\mat{I}-\lambda \mat{P})^{-1} \mat{r}.
\end{equation}




For the purposes of sensitivity analysis, we consider $R_t$ and $R_\infty$ as functions of $\mat{\pi}$, $\mat{P}$, and $\mat{r}$, since these are the quantities most likely to be estimated from data and subject to uncertainty. We assume that $\lambda$ is fixed. Clearly, $R_t$ is a polynomial in $\mat{\pi}$, $\mat{P}$, and $\mat{r}$. For $R_\infty$, we can re-write the matrix inversion as
\begin{equation}
    \label{eq:infty_reward_nice}
    R_\infty = \frac{1}{\det{(\mat{I} - \lambda \mat{P})}} \mat{\pi}^\top  \adj{(\mat{I} - \lambda \mat{P})} \mat{r},
\end{equation}
where $\det(\cdot)$ is the determinant and $\adj(\cdot)$ is the adjugate matrix operator \citep{stranglinalg}. In this form, there are three key properties of $R_\infty$ that we will use.



\begin{lemma}
    \label{lem:infy_reward_ratio}
    $R_\infty$ is a ratio of two polynomials in $\mat{\pi}$, $\mat{P}$, and $\mat{r}$. Furthermore, $\det(\mat{I} - \lambda \mat{P}) > 0$ and the adjugate $\adj(\mat{I} - \lambda \mat{P})$ has all non-negative entries. 
\end{lemma}

\proof{Proof.}
    Let $\mat{M} = \mat{I} - \lambda \mat{P}$. Each element of the adjugate matrix is the determinant of a submatrix of $\mat{M}$, and the determinant is a multilinear map, and so each element of $\adj{\mat{M}}$ is a polynomial of a submatrix of $\mat{P}$. Thus $\mat{\pi}^\top \adj{(\mat{M})} \mat{r}$ is a polynomial of $\mat{\pi}$, $\mat{P}$, and $\mat{r}$. As noted, $\det{\mat{M}}$ is a polynomial of $\mat{P}$, so $R_\infty$ is a ratio of two polynomials. Next, $\mat{M}$ has the property that all of its real eigenvalues are positive \citep[Theorem 2.3 in Chapter 6]{berman1994matrices}. Since any complex eigenvalues come in conjugate pairs, the product of all eigenvalues, which equals the determinant, is positive. Lastly, the inverse of $\mat{M}$ has all non-negative entries \citep[Theorem 2.3 in Chapter 6]{berman1994matrices}, and we already know its determinant is positive, so the adjugate of $\mat{M}$ has all non-negative entries.
    \halmos
\endproof

\begin{example}
    In a Markov chain with $n=2$ states, $\det(\mat{I} - \lambda \mat{P}) = 1 + \lambda^2 p_{11} p_{22} - \lambda^2 p_{12} p_{21} - \lambda p_{11} - \lambda p_{22}$ and $\adj(\mat{I} - \lambda \mat{P})$ is
    \begin{equation*}
    \begin{bmatrix}
    1-\lambda p_{22} & \lambda p_{12} \\
    \lambda p_{21} & 1-\lambda p_{11} 
    \end{bmatrix}.
    \end{equation*}
    Hence, $R_\infty$ is the following ratio of polynomials in $\mat{P}$, $\mat{\pi}$, and $\mat{r}$: 
    \begin{align*}
        \frac{r_1 (\lambda p_{21} \pi_2 + \pi_1 (1 -\lambda p_{22})) + r_2(\lambda p_{12} \pi_1 + \pi_2 (1 -\lambda p_{11}))}{1 + \lambda^2 p_{11} p_{22} - \lambda^2 p_{12} p_{21} - \lambda p_{11} - \lambda p_{22}}.
    \end{align*}
\end{example} 


Sensitivity analyses of Markov reward processes, and cost-effectiveness analyses in particular, are typically concerned with identifying a range of input values over which inequalities involving $R_t$ and $R_\infty$ hold. Given the forms of $R_t$ and $R_\infty$ provided above, particularly in Lemma \ref{lem:infy_reward_ratio} for $R_\infty$, typical inequalities associated with sensitivity analysis can all be written as polynomial inequalities. In the examples that follow, we focus on $R_\infty$; the application to $R_t$ is straightforward. Importantly, since the denominator of $R_\infty$ is always positive, the sign of the inequality remains unchanged as elements of $\mat{P}$ are varied. 

\paragraph{Total reward.} Given a threshold $T \in \mathbb{R}$, a policymaker may be interested in
\begin{equation}
    \mat{\pi}^\top (\mat{I}-\lambda \mat{P})^{-1} \mat{r} \geq T,
\end{equation}
which can be written as the polynomial inequality
\begin{equation}\label{eq:totalrewardineq}
    \mat{\pi}^\top \adj(\mat{I} - \lambda \mat{P}) \mat{r} - T \det(\mat{I} - \lambda \mat{P}) \geq 0.
\end{equation}

This approach can be extended to the comparison of two interventions, labeled $a$ and $b$. Determining whether the total reward associated with intervention $a$ is greater than that of $b$ 
\begin{equation}
\mat{\pi}_a^\top (\mat{I} - \lambda \mat{P}_a)^{-1} \mat{r}_a \geq \mat{\pi}_b^\top (\mat{I} - \lambda \mat{P}_b)^{-1} \mat{r}_b
\end{equation}
can be written as the polynomial inequality
\begin{equation}
\label{eq:reward_comparison_full}
    \mat{\pi}_a^\top \adj(\mat{I} - \lambda \mat{P}_a) \mat{r}_a \det(\mat{I} - \lambda \mat{P}_b) - \mat{\pi}_b^\top \adj(\mat{I} - \lambda \mat{P}_b) \mat{r}_b \det (\mat{I} - \lambda \mat{P}_a) \geq 0.
\end{equation}
Comparing more than two interventions would result in a system of polynomial inequalities. 

\paragraph{Net monetary benefit (NMB).} The net monetary benefit is defined as a constant willingness-to-pay threshold ($W$) for one unit of benefit, multiplied by the infinite horizon benefit, and then subtracting the infinite horizon cost. Let $\mat{b}$ and $\mat{c}$ be the vectors representing the one-period benefit and cost for each state. Then NMB equals
\begin{equation}
\label{eq:nmb}
    W \mat{\pi}^\top (\mat{I}-\lambda \mat{P})^{-1} \mat{b} - \mat{\pi}^\top (\mat{I} - \lambda \mat{P})^{-1} \mat{c} = \mat{\pi}^\top (\mat{I} - \lambda \mat{P})^{-1} (W\mat{b}-\mat{c}).
\end{equation}
The right-hand side expression is equivalent to the total reward when $\mat{r} = W\mat{b}-\mat{c}$. Hence, an inequality that bounds the NMB can be written similarly to \eqref{eq:totalrewardineq}. 

\paragraph{Incremental cost-effectiveness ratio (ICER).} 
The incremental cost-effectiveness ratio is defined as the ratio between the difference in the infinite horizon costs and the difference in the infinite horizon benefits of two interventions. For two interventions $a$ and $b$, the ICER is
\begin{equation}
\label{eq:icer}
\frac{\mat{\pi}_a^\top (\mat{I}-\lambda \mat{P}_a)^{-1} \mat{c}_a - \mat{\pi}_b^\top (\mat{I}-\lambda \mat{P}_b)^{-1} \mat{c}_b}{\mat{\pi}_a^\top (\mat{I}-\lambda \mat{P}_a)^{-1} \mat{b}_a - \mat{\pi}_b^\top (\mat{I}-\lambda \mat{P}_b)^{-1} \mat{b}_b}.
\end{equation}
Bounding the ICER by $T$ results in a polynominal inequality after rearranging terms. 

\begin{remark} [Death state]
    Many Markov chains used in health economic models include an absorbing ``death'' state: once the process transitions there it remain in this state with probability $1$ and reward $0$. 
    In this case, the infinite horizon expected total reward \textit{without} discounting is finite and can be written as \citep{putermanmarkov}:
\begin{equation}
\label{eq:deathstate}
\mat{\bar{\pi}}^\top (\mat{I} - \mat{Q})^{-1} \mat{\bar{r}},
\end{equation}
    where $\mat{Q}$ is a $(n-1) \times (n-1)$ matrix representing transitions between the transient states, and $\mat{\bar{r}}$ and $\mat{\bar{\pi}}$ are the subvectors corresponding to the rewards and initial distribution, respectively, on the transient states. Bounding this quantity can be similarly reformulated into a polynomial inequality. Note that the constraints on $\mat{Q}$ will be \textit{substochastic}, i.e., the row sums need only be less than or equal to 1. As well, the row sums of $\mat{Q}$ need to be strictly greater than zero to ensure that $\mat{I} - \mat{Q}$ is invertible. Multiple absorbing states can also be accommodated easily. The takeaway is that our subsequent development for discounted Markov reward processes applies to this case as well. 

\end{remark}


Beyond the polynomial inequalities above associated with the total reward, constraints on the inputs $\mat{\pi}$, $\mat{P}$, and $\mat{r}$ may contribute additional polynomial inequalities. Examples include simplex constraints on $\mat{\pi}$ and $\mat{P}$, and application-specific constraints such as non-negativity of $\mat{r}$ or $\mat{P}$ possessing the increasing failure rate property. Thus, for both finite horizon and infinite horizon, a complete sensitivity analysis can be done by examining a system of several polynomial inequalities. 


\section{Cylindrical Algebraic Decomposition}
\label{sec:cad}

Motivated by the question \textit{``for what values of $\mat{\pi}$, $\mat{P}$, and $\mat{r}$ does a given inequality hold?''}, this section studies general systems of polynominal inequalites. The focus is on those systems that possess relevant characteristics associated with Markov reward process analysis, as described in the previous section.  For solving polynomial inequalities (i.e., identifying regions or points in the domain that satisfy the inequalities), the standard representation is a \textit{cylindrical algebraic decomposition} (CAD). We provide an overview of CAD following \citet{basu} and refer the reader to that book for more details. 




A CAD of a polynomial system is a finite decomposition of $\mathbb{R}^k$ into disjoint \textit{cells} where in each cell each polynomial in the system is sign-invariant\endnote{There can be many different CADs to represent the same system, e.g., by splitting cells unnecessarily.}. Having such a decomposition allows us to easily test consistency of the system at a given point in the domain and describe the regions over which the system is consistent.  Next, we formally define a cell. 

 
\begin{definition}[Cell]
    A cell is defined recursively.
    \begin{enumerate}
        \item In $\mathbb{R}^1$, a cell is either an open interval or a point.
        \item Let $k \ge 1$ and $C$ be a cell in $\mathbb{R}^k$. In $\mathbb{R}^{k+1}$, a cell is either of the form $\{(x,y) \in \mathbb{R}^{k+1} \;|\; x \in C, f(x) < y < g(x)\}$ or $\{(x,y) \in \mathbb{R}^{k+1} \;|\; x \in C,  y = f(x)\}$, where $f$ and $g$ are either algebraic functions\endnote{An algebraic function is defined as being the zero of some polynomial. It is a strictly larger class of functions than the polynomials. For example, $\sqrt{x}$ is not a polynomial, but it is algebraic, as it is a solution to $y^2 - x = 0$. On the other hand, $\sin{x}$, which is certainly not a polynomial, is also not algebraic.} or $\pm \infty$, with $f(x) < g(x)$ for all $x \in C$. 
    \end{enumerate}
\end{definition}

Let $\mathbb{R}[\{x_i\}_{i=1}^k]$ be the ring of polynomials in $x_1, \ldots, x_k$. If we have a finite set of polynomials $F \subset \mathbb{R}[\{x_i\}_{i=1}^k]$, we call a CAD \textit{adapted to} $F$ if each $f \in F$ is sign-invariant ($0$, $+$, or $-$) over each cell. An important theorem from algebraic geometry is that for every finite set of polynomials $F$, there exists a CAD adapted to $F$ (see \citet{lojasiewicz1965ensembles} or \citet[Theorem 5.6]{basu}).

Next, we consider inequalities involving elements of $F$. Let $f \in \mathbb{R}[\{x_i\}_{i=1}^k]$. An \textit{$f$-atom} is an expression $f \bowtie 0$, where $\bowtie \in \{=, >, <, \geq, \leq\}$. A \textit{semialgebraic set} is defined by Boolean combinations of $f$-atoms and is closed under union, intersection and complement, . 
A corollary of the above theorem is that any semialgebraic set of $f$-atoms can be represented equivalently by a subset of the cells of the CAD adapted to $F$. 

The recursive nature of the cells of a CAD make it intuitive to represent as a tree.

\begin{example} 
    \label{ex:sphere_cad}
    Consider the unit sphere $x^2+y^2+z^2-1=0$. Its simplest CAD is:.
    \begin{figure}[H]
        \centering
        {\scriptsize
\begin{forest}
[{}, phantom
[{$x=-1$}
    [{$y=0$}
      [{$z=0$}]
    ]
]
[{$-1<x<1$}
    [{$y=-\sqrt{1-x^2}$}
        [{$z=0$}]
    ]
    [{$-\sqrt{1-x^2}<y<\sqrt{1-x^2}$}
        [{$z=-\sqrt{1-x^2-y^2}$}]
        [{$z=\sqrt{1-x^2-y^2}$}]
    ]
    [{$y=\sqrt{1-x^2}$}
        [{$z=0$}]
    ]
]
[{$x=1$}
    [{$y=0$}
        [{$z=0$}]
    ]
]
]
\end{forest}
}
    \end{figure}
    The CAD can be constructed using straightforward geometric reasoning. We start with one-dimensional cells, so we project the sphere to the unit disc $x^2+y^2 \leq 1$ in $\mathbb{R}^2$, and then again to the interval $-1 \le x \le 1$ in $\mathbb{R}^1$. Form cells $-1$, $(-1,1)$, and 1. Given the $x$ values in each cell, determine the feasible $y$ values based on the unit disc. Repeat for the $z$ values based on the unit sphere.
    
    A cell in $\mathbb{R}^3$ is the conjunction of all the nodes in a path from $x$ to $z$. An example is $x=-1 \land y=0 \land z=0$. By taking the disjunction of all cells, we get the full decomposition of the unit sphere.
\end{example}

This simple example illustrates the general principles for how to construct a CAD for any polynomial system in $\mathbb{R}^k$: 1) Repeatedly project the polynomials to $\mathbb{R}^1$; 2) Create cells by identifying appropriate points (i.e., roots) and intervals; 3) Over each cell in $\mathbb{R}^1$, construct a cylinder in $\mathbb{R}^2$, and create cells in $\mathbb{R}^2$ by computing the roots and intervals of each projection in $\mathbb{R}^2$. 4) Repeat until we reach $\mathbb{R}^k$. 

While the above approach is intuitive, it required the development of special projection operators \citep{collins1974quantifier,collins1976quantifier} to allow systematic CAD construction for arbitrary systems. Next, we describe the algorithm for generating a CAD of a system.



\subsection{The CAD Algorithm}

The algorithm of \citet{collins1974quantifier} is given below. We refer the reader to Chapters 5 and 11 of \citet{basu} for a comprehensive description of the algorithm. Here, we break up CAD into two algorithms, which we call the \textit{decision phase} and \textit{solution formula phase}. The pseudocode of these algorithms is given in Algorithm \ref{alg:caddecision} and Algorithm \ref{alg:cadsoln}, respectively. The decision phase of CAD determines whether a polynomial system is consistent, whereas the solution formula phase constructs the full CAD tree as in Example \ref{ex:sphere_cad}.

\begin{algorithm}
\caption{Cylindrical algebraic decomposition: decision phase}
\label{alg:caddecision}
\begin{algorithmic}[1]
\Require A set of polynomials $F_k \subset \mathbb{R}[\{x_i\}_{i=1}^k]$
\Ensure Sample points in each cell of the CAD adapted to $F_k$

\Statex \textbf{\# Step 1: Projection}
\For{$i = k, k-1, \dots, 2$}
    \State From $F_{i}$, construct the set of projection factors $F_{i-1} \subset \mathbb{R}[\{x_j\}_{j=1}^{i-1}]$ by eliminating $x_i$
\EndFor

\Statex \textbf{\# Step 2: Base Case}
\State $r \gets$ Ordered list of roots of all $f \in F_1$, labelled $r_1, r_2, \ldots$
\State Construct cells: $\{r_j\}_{j=1}^{|r|} \cup \{(r_j, r_{j+1})\}_{j=1}^{|r|-1} \cup \{(-\infty, \min r), (\max r, \infty)\}$
\Statex \textit{\# Calculate sample points}
\For{each cell}
    \If{cell is a point}
        \State Sample point is the point itself
    \ElsIf{cell is an interval with finite endpoints}
        \State Sample point is average of endpoints
    \ElsIf{cell is $(-\infty, \min r)$}
        \State Sample point is $\min r - 1$
    \ElsIf{cell is $(\max r, \infty)$}
        \State Sample point is $\max r + 1$
    \EndIf
\EndFor

\Statex \textbf{\# Step 3: Lifting}
\State $\mathcal{C} \gets$ sample points of cells in $\mathbb{R}^1$
\For{$i = 2, 3, \dots, k$}
    \State $\mathcal{C'} \gets \emptyset$
    \For{each cell $C \in \mathcal{C}$}
        \State Evaluate the polynomials in $F_i$ at the sample point of $C$
        \State Find the roots of these polynomials and construct cells and sample points in $\mathbb{R}^i$ as in the base case
        \State Add these sample points to $\mathcal{C'}$ 
    \EndFor
    \State $\mathcal{C} \gets \mathcal{C'}$
\EndFor

\State \Return $\mathcal{C}$

\end{algorithmic}
\end{algorithm}

The input of the decision phase is a set $F_k \subset \mathbb{R}[\{x_i\}_{i=1}^k]$ of polynomials in $k$ variables. The algorithm has three main steps: projection, base case, and lifting. In the projection step, the set $F_k$ is iteratively projected down to lower dimensions, forming sets of polynomials that are called \textit{projection factors} in $\mathbb{R}^{k-1}, \mathbb{R}^{k-2}, \cdots \mathbb{R}^1$, denoted $F_{k-1}, F_{k-2}, \cdots F_1$, respectively. In the base case step, we form cells in $\mathbb{R}^1$ using $F_1$ and then store a single sample point for each cell. Finally, in the lifting step, the CAD is lifted iteratively back up to $\mathbb{R}_k$. The final result is a set of points in $\mathbb{R}^k$ where each point represents a cell over which each $f \in F_k$ is sign-invariant. Therefore, the truth of a system of polynomials where the atoms are composed of $f \in F_k$ can be determined by evaluation at each test point, due to the sign-invariance property. 
Thus, Algorithm \ref{alg:caddecision} allows us to determine if \textit{any} solution exists, which is sufficient for many applications, e.g., SMT. 

If we want an algebraic description of \textit{all} points that satisfy the system, like in the case of a sensitivity analysis, then we require Algorithm \ref{alg:cadsoln}. This algorithm produces the desired result for \textit{projection-definable} systems \citep{brownthesis}. A system is projection-definable if no two cells that have different truth values share the same sign for any projection factor. In this case, the solution formula can be constructed purely using the signs of the projection factors. 
Algorithm \ref{alg:cadsoln} requires information from Algorithm \ref{alg:caddecision}, namely the set of all projection factors and cells (specifically, their sample points), and relies on the sign-invariance of projection factors over the cells in their respective dimension. Algorithm \ref{alg:cadsoln} is sufficient, because as we will show in Lemma \ref{lem:xg_projdef}, all polynomial systems that arise in our sensitivity analysis application are projection-definable. Systems that are not projection-definable are outside the scope of this paper. The interested reader is referred to \citet{brownthesis}. 


\begin{algorithm}
\caption{Cylindrical algebraic decomposition: solution formula phase}
\label{alg:cadsoln}
\begin{algorithmic}[1]
\Require A projection-definable polynomial system, and all projection factors $F_k \cup F_{k-1} \cup F_{k-2} \cup \cdots \cup F_1$ and the set of sample points $\mathcal{C}$ from Algorithm \ref{alg:caddecision}
\Ensure A CAD-based formula representing all the solutions to the polynomial system

\State \textit{cell\_formulas} $\gets$ empty list
\For{$C \in \mathcal{C}$ where the polynomial system holds over $C$}
    \State \textit{formulas} $\gets$ empty list
    \For{$f \in F_k \cup F_{k-1} \cup F_{k-2} \cdots \cup F_1$}
        \If{$f$ is negative at $C$}
            \State Append atom $f < 0$ to \textit{formulas}
        \ElsIf{$f$ is zero at $C$}
            \State Append atom $f = 0$ to \textit{formulas}
        \ElsIf{$f$ is positive at $C$}
            \State Append atom $f > 0$ to \textit{formulas}
        \EndIf
    \EndFor
    \State Append $\bigwedge$ \textit{formulas} to \textit{cell\_formulas}
\EndFor
\State \Return $\bigvee$ \textit{cell\_formulas}

\end{algorithmic}
\end{algorithm}

The crucial step in Algorithm \ref{alg:caddecision} is the projection step (line 2), which is now the most well-studied part of the algorithm. The first valid projection operator was given by \citet{collins1974quantifier} and later simplified by \citet{hong1990improvement}, which gives smaller projection factor sets while still maintaining soundness and completeness. 
In our implementation, we use the \citet{hong1990improvement} projection operator, which is defined in Appendix \ref{sec:hong}. In Algorithm \ref{alg:cadsoln}, the conjunctions in line 13 can be visualized in the tree in Example \ref{ex:sphere_cad} as the paths from the top of the tree to the leafs, while the disjunctions in line 15 are the set of all paths, which comprise the entire tree.

We now have all the information necessary to rigorously construct the CAD of an arbitrary system. As an example, we revisit the 3D sphere CAD and construct it per Algorithms \ref{alg:caddecision} and \ref{alg:cadsoln}. 

\begin{example}
\label{ex:sphere_cad_execution}
    For the unit sphere, our system is $x^2+y^2+z^2-1=0$, so $F_3 = \{x^2+y^2+z^2-1\}$. We eliminate $z$ and then $y$.

    \emph{Algorithm \ref{alg:caddecision}:}

    \begin{hitemize}
        \item Projection (see Appendix \ref{sec:hong} for a detailed execution of the Hong projection operator)
        \begin{hitemize}
            \item Eliminate $z$: $F_2 = \{x^2+y^2-1, -4x^2-4y^2+4\}$
            \item Eliminate $y$: $F_1 = \{-256+256x^2, -4x^2+4, x^2-1, x^4-2x^2+1\}$
        \end{hitemize}

        \item Base case. The roots of $F_1$ are $\{-1, 1\}$. So the cells are $\{(-\infty, -1), -1, (-1, 1), 1, (1, \infty)\}$. We choose the corresponding sample points $\{-2, -1, 0, 1, 2\}$.

        \item Lifting
        \begin{hitemize}
            \item To $\mathbb{R}^2$:
            \begin{hitemize}
                \item Over sample point $-2$, we have the polynomials $\{3+y^2, 12-4y^2\}$, with roots $\{-\sqrt{3}, +\sqrt{3}\}$, hence cells $\{(-\infty, -\sqrt{3}), -\sqrt{3}, (-\sqrt{3}, \sqrt{3}), \sqrt{3}, (\sqrt{3}, \infty)\}$, with sample points $\{-\sqrt{3}-1, -\sqrt{3}, 0, \sqrt{3}, \sqrt{3}+1\}$
                \item Over sample point $-1$, we have the polynomials $\{y^2, -4y^2\}$, with roots $\{0\}$, hence cells $\{(-\infty, 0), 0, (0, \infty)\}$, with sample points $\{-1, 0, 1\}$
                \item Over sample point $0$, we have the polynomials $\{y^2-1, -4y^2+4\}$, with roots $\{-1, 1\}$, hence cells $\{(-\infty, -1), -1, (-1,1), 1, (1, \infty)\}$, with sample points $\{-2, -1, 0, 1, 2\}$
                \item Over sample points $1$ and $2$, the analysis is the same as with sample point $-1$ and $-2$, respectively, because $x$ always appears as $x^2$ in $F_1$
            \end{hitemize}

            \item To $\mathbb{R}^3$: We plug in each $(x,y)$ sample point into $x^2+y^2+z^2-1$ and construct cells and sample points for $z$; we omit the details for brevity. There are $41$ cells in total.
        \end{hitemize}
        \item Truth of the system: Evaluate $x^2+y^2+z^2-1$ at each sample point and check if it equals zero. The system holds over the following $6$ cells, represented by their $(x,y,z)$ sample points: $\{(-1,0,0), (0,-1,0), (0,0,-1), (0,0,1), (0,1,0), (1,0,0)\}$.
    \end{hitemize}

    \emph{Algorithm \ref{alg:cadsoln}:}

    The set of projection factors is $F_1 \cup F_2 \cup F_3$. We evaluate their signs over each cell. We can ignore $\{-4x^2-4y^2+4, -4x^2+4, -256+256x^2\}$ as they are multiples of other projection factors. We can also ignore $x^4-2x^2+1$, because it equals $(x^2-1)^2$, so it is always non-negative, and zero if and only if $-1+x^2$ is zero, so it is also redundant to list in the solution formula.\endnote{Note that these `algebraic reasoning' arguments like ignoring multiples, solving inequalities of univariate polynomials, etc., are not strictly part of CAD -- they should be handled at the level of a (good) computer algebra system.} Thus, the three projection factors that remain are $\{x^2-1, x^2+y^2-1, x^2+y^2+z^2-1\}$. This system is projection-definable, as no pair of cells with differing truth values share the same signs of the projection factors, which can be manually checked. 
    
    Next, we show an example of constructing the solution formula for a given a cell. Over the cell with sample point $(0,0,1)$ applied to the three remaining projection factors we get the atoms $\{x^2-1<0, x^2+y^2-1<0, x^2+y^2+z^2-1=0\}$, and hence the formula $-1<x<1 \land -\sqrt{1-x^2}<y<\sqrt{1-x^2} \land z = \sqrt{1-x^2-y^2}$. Following a similar approach for the other cells, we recover the full CAD of the sphere in Example \ref{ex:sphere_cad}.
\end{example}

Although we have a general approach to construct the CAD of an arbitrary system, the key challenge is computational complexity. Due to the complicated ways polynomials can interact during the projection phase, the number of projection factors and the size of the CAD can quickly grow large. CAD has a complexity that depends polynomially on the degrees and number of polynomials, but is \textit{doubly exponential} in the number of variables \citet{england2015improving}. Indeed, it is possible to construct examples where this doubly exponential complexity is attained \citep{basu}. This presents a significant challenge in computing CADs for even more than a few variables. As well, the order of which variables to eliminate during the projection step can significantly affect both the runtime and the space needed to store the CAD. 


\subsection{A Special Class of Polynomial System}\label{sec:specialclass}

In the remainder of this section we examine a specific class of polynomial systems that is motivated directly by our Markov reward process context. We demonstrate that due to the special structure of this system, its CAD can be generated much more efficiently than the general case. We also present a specialized version of the general CAD algorithm, tailored to this polynomial system class. 

We consider a system of polynomials that have two different ``types'' of variables: $x$-type variables and $\alpha$-type variables. We have $\eta$ of the $x$-type variables, so that we index them $\{x_i\}_{i=1}^\eta$. The $\alpha$-type variables are doubly indexed and are arranged as follows: there are $\phi$ simplices, and in the $i$th simplex, there are $\tau_i$ of the $\alpha$-type variables. In other words, we index the $\alpha$-type variables as $\{\alpha_{i,j}\}_{i \in [\phi], j \in [\tau_i])}$, where we define the notation $[m] := \{1, \ldots, m\}$ for a positive integer $m$. By definition, each variable $\alpha_{i,j}$ is present in only a single simplex. Thus, there are a total of $\eta+\sum_{i=1}^\phi \tau_i$ variables. Let $\mathbb{R}[\{x_i\}_{i=1}^k,\{\alpha_{i,j}\}_{i \in [\phi], j \in [\tau_i]}]$ represent the ring of polynomials of these variables, and let $f^*$ be a polynomial in this ring. We will call the atom $f^* \geq 0$, the \textit{defining inequality} of the system. The connection to our sensitivity analysis motivation is that $f^* \ge 0$ is the condition the policymaker wants to test (e.g., NMB greater than a threshold), the $x$-type variables represents rewards, costs and benefits, and the $\alpha$-type variables represent probabilities, i.e., $\mat{\pi}$ and the rows of $\mat{P}$.


Our polynomial system of interest is

\begin{equation}
\label{eq:general_form}
\begin{aligned}
    & f^* \geq 0  \\
    & x_i \geq 0 \quad i \in [\eta]  \\ 
    & \sum_{j=1}^{\tau_i} \alpha_{i,j}=1 \quad i \in [\phi] \\ 
    & 0 \leq \alpha_{i,j} \leq 1 \quad i \in [\phi], j \in [\tau_i]
\end{aligned}
\tag{M}
\end{equation}

We focus on the case where $x_i$ is sign-constrained both because of  practical motivations (e.g., rewards and costs associated with real-world problems tend to be sign-constrained) and because it leads to slightly more complex CAD trees (due to the presence of the endpoint 0) allowing us to showcase the full range of the CAD algorithm. The following development fully applies to the case where $x_i$ is free, except that there will be fewer cells.

Despite the general challenges described at the end of the previous subsection, we will show that we can construct a CAD for system \eqref{eq:general_form} far more efficiently than general systems. This is made possible due to the following:

\begin{enumerate}
    \item For a given $i$, the $\alpha_{i,j}$ variables form a unit simplex.
    \item For a given $i,j$, the $\alpha_{i,j}$ variable only appears in a single simplex constraint.
    \item The $x$-type variables only appear in $f^* \ge 0$.
\end{enumerate}

Our approach takes advantage of this structure as follows. First, we construct a CAD of each simplex (which can be parallelized), and then conjunct them together. Then, relying on a special property associated with the projection factors of $f^*$ that we describe later (``simplex-extensibility''), we can efficiently lift to include the $x$-type variables, which only show up in the defining inequality, to construct the full CAD.

\subsubsection{CAD of simplices.}

Since all variables are restricted to be between $0$ and $1$, the construction of a simplex CAD is straightforward.

\begin{theorem}
\label{thm:simplex_cad}
    Consider a unit simplex $\{\alpha_j, j \in [\tau] \,|\, \sum_{j=1}^\tau\alpha_j=1, 0 \le \alpha_j \le 1, j \in [\tau]\}$. Its CAD is given by the following tree:
    \begin{enumerate}
        \item The first level has three cells: $\alpha_1=0$, $0<\alpha_1<1$, $\alpha_1=1$.
        \item For $2 \leq i \leq \tau - 1$, for a given cell in level $i-1$, the children cells in level $i$ are $0$, $(0,1-\sum_{j=1}^{i-1} \alpha_j)$ and $1-\sum_{j=1}^{i-1} \alpha_j$.
        \item The last level ($i=\tau$) has a single child cell for each of the $\tau-1$ level cells: $\alpha_\tau = 1-\sum_{j=1}^{\tau-1} \alpha_j$.
    \end{enumerate}
    The tree also propagates down any cells defined by an equality on the preceding variables.
\end{theorem}

\begin{remark}
    By propagating any equality constraints active on the current cell, we may obtain a degenerate interval (i.e., a point) in which case there is only one new cell. This is the case if, for example, (taking a simplex with four variables) $0 < \alpha_1 < 1 \land \alpha_2 = 1-\alpha_1$. Then, surely $\alpha_3=\alpha_4=0$.
\end{remark}

It is easy to extend this simplex construction to the case where $\sum_{i=1}^\tau \alpha_i = \kappa \leq 1$, as the defining representation of the simplex is now $\alpha_1 \in [0,1]$, $\alpha_2 \in [0, \kappa - \alpha_1]$, $\alpha_3 \in [0, \kappa - \alpha_1 - \alpha_2]$, and so on, until $\alpha_\tau = \kappa - \sum_{i=1}^{\tau-1} \alpha_i$. Also, we can extend this construction to the inequality case, where $\sum_{i=1}^\tau \alpha_i \leq \kappa \leq 1$, by replacing the node at the final level with $\alpha_\tau \leq \kappa - \sum_{i=1}^{\tau - 1} \alpha_i$.

Since each $\alpha$-type variable only occurs in a single simplex, it is trivial to create a CAD for a set of simplices, by ``gluing" them together.

\begin{corollary}
\label{cor:gluing_simplex_cad}
    Consider $\phi$ unit simplices $\{\alpha_{i,j}, i \in [\phi], j \in [\tau_i] \,|\, \sum_{j=1}^{\tau_i} =1, 0 \le \alpha_{i,j} \le 1, i \in [\phi], j \in [\tau_i]\}$.
    To construct a cell in the CAD of the conjunction of these simplices, choose a single cell from each of the individual simplices' CADs, and conjunct them. The full CAD is the disjunction of all such cells.
\end{corollary}

\begin{example}
    Consider two simplices $\alpha_{1,1} + \alpha_{1,2} + \alpha_{1,3} = 1$ and $\alpha_{2,1} + \alpha_{2,2} + \alpha_{2,3} = 1$. The CAD of their conjunction is:

    \begin{figure}[H]
        \centering
        {\scriptsize
\begin{forest}
for tree={baseline}
[{}, phantom
[{$\alpha_{1,1} = 0$}
	[{$\alpha_{1,2} = 0$}
		  [{$\alpha_{1,3} = 1$}, tikz={\node[draw,fit=(first) (first-1) (second) (second-1) (third) (third-1), label=below:S] {};}
                [{$\alpha_{2,1} = 0$}, l=2cm, name=first [{}, l=1cm, edge=dashed, name=first-1]]
                [{$0 < \alpha_{2,1} < 1$}, l=2cm, name=second [{}, l=1cm, edge=dashed, name=second-1]]
                [{$\alpha_{2,1} = 1$}, l=2cm, name=third [{}, l=1cm, edge=dashed, name=third-1]]
            ]
	]
	[{$0 < \alpha_{1,2} < 1$}
		  [{$\alpha_{1,3} = 1 - \alpha_{1,2}$} [S]]
	]
	[{$\alpha_{1,2} = 1$}
		  [{$\alpha_{1,3} = 0$} [S]]
	]
]
[{$0 < \alpha_{1,1} < 1$}
	[{$\alpha_{1,2} = 0$}
		  [{$\alpha_{1,3} = 1 - \alpha_{1,1}$} [S]]
	]
	[{$0 < \alpha_{1,2} < 1 - \alpha_{1,1}$}
		  [{$\alpha_{1,3} = 1 - \alpha_{1,1} - \alpha_{1,2}$} [S]]
	]
	[{$\alpha_{1,2} = 1 - \alpha_{1,1}$}
		  [{$\alpha_{1,3} = 0$} [S]]
	]
]
[{$\alpha_{1,1} = 1$}
	[{$\alpha_{1,2} = 0$}
		  [{$\alpha_{1,3} = 0$} [S]]
	]
]
]
\end{forest}
}
        \vspace{-1em}
    \end{figure}

    \noindent where $S$ represents the CAD of the second simplex, $\alpha_{2,1} + \alpha_{2,2} + \alpha_{2,3} = 1$, shown partially on the leftmost branch, which has identical structure to the CAD of the first simplex.
\end{example}

The intuition of the preceding corollary is that because each simplex is independent, we can simply copy and conjunct them together. For example, if we have two simplices, take each cell in the first simplex, and attach a copy of the second simplex. This leads to a bound on the size of the CAD of simplices.

\begin{corollary}
    \label{cor:simplex_cad_numcells}
    The number of cells in a CAD of the conjunction of $\phi$ simplices, where the $i$th simplex has $\tau_i$ variables, is $O(3^{\sum_{i=1}^\phi \tau_i}) = O(3^{\phi  \max_i \tau_i})$.
\end{corollary}

In general, the size of a CAD is doubly exponential in the number of variables, whereas for special case of simplices it is singly exponential: the savings are due to the geometry of the simplex and the assumption of their disjointness.

\subsubsection{Simplex-extensibility: Lifting the simplex CAD to include $f^*$.}

The previous section showed that it was easy to construct a CAD for a set of simplices. Now, we wish to lift the CAD to include $f^*$. The method of incremental CAD, which we discussed in the literature review, applies here. In the general case, it requires significant computation, due to the many possibilities of how CADs can interact. Therefore, we focus on characterizing conditions on $f^*$ such that extending the simplex CAD is easy.

Suppose we repeatedly apply a projection operator to the set $F$ of polynomials that make up our system \eqref{eq:general_form} to eliminate all $x$-type variables. Then, we are left with a projection factor set $F'$, where each $f \in F'$ is a polynomial in the $\alpha$-type variables. If the CAD of the simplex constraints is \textit{also} the CAD of $F'$, then we can simply stop the projection phase here and begin the lifting procedure from the simplex CAD to construct a CAD for the full system. The value of such an approach is the avoidance of calculating subsequent projection factor sets, which can blow up quickly. Thus, we wish to characterize the polynomials $f^*$ such that this is true.

To do so, we introduce several technical lemmas in Appendix \ref{sec:some_technical_lemmas} that help us to apply the Hong projection operator to instances of system \eqref{eq:general_form}. Below, we provide a motivating example that will demonstrate the situation we discussed above.

\begin{example}
    \label{ex:projs}
    Consider the system $\alpha_1 x_1 + \alpha_2 x_2 -1 \geq 0$ and $\alpha_1+\alpha_2=1$. The set of polynomials whose atoms form the system are $\{\alpha_1 x_1 + \alpha_2 x_2 -1, \alpha_1+\alpha_2-1\}$. Iteratively applying the Hong projection:
    
    \begin{itemize}
        \item Eliminate $x_1$: $\{\alpha_1, \alpha_1 x_2-1, \alpha_1+\alpha_2-1\}$
        \item Eliminate $x_2$: $\{\alpha_1,\alpha_2,\alpha_1+\alpha_2-1\}$
    \end{itemize}

    Now, observe that the final projection factor set has the same CAD as the CAD of the simplex associated with $\alpha_1 + \alpha_2 = 1$, so we can lift from it. Specifically, observe that each of the elements of the projection factor set are \textit{sign-invariant} in each of the cells. 
\end{example}

The sign-invariance property observed in the previous example is crucial to our development. We call this property \textit{simplex-extensible}.

\begin{definition}[Simplex-extensible]
    An instance of the system \eqref{eq:general_form} is \textit{simplex-extensible} if the set of projection factors after eliminating all $x$-type variables is sign-invariant over each cell in the CAD of the conjunction of the respective simplex constraints.
\end{definition}

\begin{example}
    \label{ex:not_simplex_ext}
    To demonstrate an instance that is not simplex-extensible, take $f^* = x_1(\alpha_2^2-\alpha_1^2)+x_2\alpha_2-1$, with the simplex constraint $\alpha_1+\alpha_2-1=0$. If we calculate the projection factor set after eliminating $x_1$ and $x_2$, we obtain $\{-1+\alpha_2, \alpha_2^2, \alpha_1+\alpha_2-1, \alpha_2^2-\alpha_1^2\}$ (we omit the details for brevity). The last polynomial, $\alpha_2^2 - \alpha_1^2$, is not sign-invariant on the simplex CAD. For example, take $(\alpha_1,\alpha_2)=(0.2,0.8)$, which satisfies the simplex: at this point, the polynomial is positive. Now take $(\alpha_1,\alpha_2)=(0.8,0.2)$, which also satisfies the simplex and indeed also lies in the same cell in the simplex CAD (namely the cell $0< \alpha_1 < 1 \land \alpha_2=1-\alpha_1$): at this point, the polynomial is negative. Therefore, it is not sign-invariant within a cell of the simplex CAD.
\end{example}

Characterizing the full class of simplex-extensible systems is difficult, as by definition it requires the repeated application of the Hong projector until all $x$-type variables are eliminated. However, using the series of technical lemmas presented in Appendix \ref{sec:some_technical_lemmas}, we can derive several special cases where the characterization is much easier. An important class of polynomials is of the form $f^* = g_0 + \sum_{i=1}^\eta x_i g_i$, where each function $\{g_i\}_{i=0}^\eta$ is a function of the $\alpha$-type variables, for which we can easily check its simplex-extensibility. 

\begin{theorem}
    \label{thm:simplex_ext_xg}
    Let $f^* = g_0 + \sum_{i=1}^\eta x_i g_i$, where each function $\{g_i\}_{i=0}^\eta$ is a polynomial of the $\alpha$-type variables. Then, this system is simplex-extensible if and only if each of $\{g_i\}_{i=0}^\eta$ is sign-invariant over the CAD of the simplex constraints.
\end{theorem}


Theorem \ref{thm:simplex_ext_xg} captures a general class of $f^*$, in particular covering the relevant functions for our Markov reward process context. More general classes of functions that exhibit simplex-extensibility are provided in Section \ref{sec:extensions_fstar}. 

In Example \ref{ex:projs}, a system with $f^* = x_1 \alpha_1 + x_2 \alpha_2 - 1$ is simplex-extensible by Theorem \ref{thm:simplex_ext_xg}. Another example is that an instance with $f^* = x_1 \alpha_1^2 + x_2 \alpha_2^2 -1$ is simplex-extensible, because $\alpha_1^2$ and $\alpha_2^2$ are strictly positive except when $\alpha_1=0$ (or likewise for $\alpha_2=0$), which constitutes its own cell in the simplex CAD. In Example \ref{ex:not_simplex_ext}, since $f^*$ doesn't satisfy the property in Theorem \ref{thm:simplex_ext_xg}, the system is not simplex-extensible. 

Thus, with a suitably chosen $f^*$, system $\eqref{eq:general_form}$ is simplex-extensible, which means we can stop projecting once we have eliminated all $x$-type variables and begin the lifting phase from the simplex CAD. Namely, we take sample points from the simplex CAD cells, plug them into the projection factors, calculate roots, and then lift on the next $x$-type variable. This now gives us a set of sample points over which we can evaluate the consistency of the system. We discuss the construction of the solution formula next.

\subsubsection{Solution formula construction.}
\label{sec:soln_construction}

Recall that we can use Algorithm \ref{alg:cadsoln} to generate the solution formula for systems that are projection-definable. Next, we show that our system of interest is projection-definable for a suitably chosen $f^*$. 


\begin{lemma}
    \label{lem:xg_projdef}
    Assume we have an instance of system \eqref{eq:general_form} that is simplex-extensible, with $f^* = g_0 + \sum_{i=1}^\eta x_i g_i$, where each function $\{g_i\}_{i=0}^\eta$ is a polynomial of the $\alpha$-type variables. Then, this system is projection-definable.
\end{lemma}

Given Lemma \ref{lem:xg_projdef}, we can specialize Algorithm \ref{alg:cadsoln} to system \eqref{eq:general_form}. We can also specialize Algorithm \ref{alg:caddecision} with knowledge of the specific projection factors associated with the simplex constraints and the structure of $f^*$. The following algorithm, Algorithm \ref{alg:cadsoln_xg}, combines the specialized versions of Algorithms \ref{alg:caddecision} and \ref{alg:cadsoln}, outputting the solution formula for any instance of \eqref{eq:general_form} with $f^* = g_0 + \sum_{i=1}^\eta x_i g_i$.

\begin{algorithm}
\caption{Constructing the CAD solution formula of an instance of system \eqref{eq:general_form} with $f^* = g_0 + \sum_{i=1}^\eta x_i g_i$}
\label{alg:cadsoln_xg}
\begin{algorithmic}[1]
\Require A simplex-extensible instance of system \eqref{eq:general_form} with $f^* = g_0 + \sum_{i=1}^\eta x_i g_i$
\Ensure A CAD solution formula

\For{each simplex constraint in the instance of system \eqref{eq:general_form}}
    \State Construct the CAD according to Theorem \ref{thm:simplex_cad}, storing sample points
\EndFor

\State \textit{cells} $\gets$ conjunction of the simplex CADs according to Corollary \ref{cor:gluing_simplex_cad}

\For{$x_i$ in $\{x_i\}_{i=1}^\eta$}
    \State \textit{cells\_new} $\gets$ \textit{cells}
    \For{each \textit{cell} in \textit{cells}}
        \State $sample$ $\gets$ sample point of \textit{cell}
        \If{$g_i = 0$, evaluated at \textit{sample}}
            \State Add the cell $x_i \geq 0$ as a child of \textit{cell} with a sample point, in \textit{cells\_new}
        \Else
            \State $r_i \gets - (g_0 - \sum_{j=1}^{i-1} x_j g_j) / g_i$ evaluated at \textit{sample}
            \If{$r_i \leq 0$}
                \State Add the cell $x_i \geq 0$ as a child of \textit{cell}, with a sample point, to \textit{cells\_new}
            \Else
                \State Add cells $x_i = 0, 0 < x_i < r_i, x_i = r_i,$ and $x_i > r_i$ as children of \textit{cell}, with their sample points, in \textit{cells\_new}
            \EndIf
        \EndIf
    \EndFor
    \State \textit{cells} $\gets$ \textit{cells\_new}
\EndFor

\For{each \textit{cell} in \textit{cells}}
    \State Evaluate the system using its sample point 
    \If{the instance of system \eqref{eq:general_form} does not hold}
        \State Delete \textit{cell}
    \EndIf
\EndFor

\State \Return \textit{cells}
\end{algorithmic}
\end{algorithm}

Line 12 of Algorithm \ref{alg:cadsoln_xg} is the key step that leverages the structure of $f^*$ to construct the CAD efficiently. At the level of a particular $x_i$, there is only one projection factor, which has a unique root (see proof of Theorem \ref{thm:simplex_ext_xg}):
\begin{equation*}
    r_i = \frac{-g_0 - \sum_{j=1}^{i-1} x_j g_j}{g_i}.
\end{equation*} 
The result is that there will be no more than four new cells created at each level of the $x$-type variables (line 16). This observation directly leads to a bound on the size of the CAD tree, and the result that the tree is smaller than the size of a CAD of a general polynomial system with the same number of variables and constraints.

\begin{corollary}
    \label{cor:full_cad_size}
    Assume we have an instance of system \eqref{eq:general_form} that is simplex-extensible, with $f^* = g_0 + \sum_{i=1}^\eta x_i g_i$, where each function $\{g_i\}_{i=0}^\eta$ is a polynomial of the $\alpha$-type variables. Then, the number of cells in its CAD is $O(3^{\sum_{i=1}^\phi \tau_i} \times 4^\eta)$.
\end{corollary}

\begin{corollary}
    \label{cor:cad_size_comparison}
    Let the number of cells of an instance of system \eqref{eq:general_form}, under the assumptions of Corollary \ref{cor:full_cad_size}, be $N_M$, and the number of cells in the CAD of a general polynomial system with the same number of variables and constraints be $N$. Then, $N_M = o(N)$.
\end{corollary}

More concretely, while the size of a general CAD is doubly exponential in size, the size of this CAD is only \textit{singly exponential}. While this is a significant reduction theoretically, we may still consider this complexity to be too high for practical purposes. However, any tree data structure with a constant number of children for each node will have a total number of nodes that is singly exponential in the number of levels of the tree. Lastly, we show that even deciding feasibility of our system (let alone enumerating the full tree), even in the simplex-extensible case, is NP-hard. The proof is via reduction from 3-SAT. 

\begin{theorem}
    \label{thm:nphard}
    Consider the decision problem of deciding the feasibility of an  instance of system \eqref{eq:general_form} that is simplex-extensible, with $f^* = g_0 + \sum_{i=1}^\eta x_i g_i$, where each function $\{g_i\}_{i=0}^\eta$ is a polynomial of the $\alpha$-type variables. This problem is NP-hard. Hence, deciding the feasibility of a general instance of system \eqref{eq:general_form} is NP-hard.
\end{theorem}


\begin{remark}[Unsigned $x$-type variables]
    Algorithm \ref{alg:cadsoln_xg} applies to the case where the $x$-type variables are free, with slight modifications to lines 10, 14, and 16. Theorem \ref{thm:simplex_ext_xg} and Lemma \ref{lem:xg_projdef} also hold because the sign of $x$ does not affect the sign-invariance of the $g_i$ functions, since the latter are functions of only the $\alpha$-type variables.
\end{remark}


The following example illustrates a CAD tree built using Algorithm \ref{alg:cadsoln_xg}. 


\begin{example}
    Let $f^* = \alpha_1 x_1 + \alpha_2 x_2 + \alpha_3 x_3 -1$, with $f^* \geq 0$, $\alpha_1+\alpha_2+\alpha_3=1$, $0 \leq \alpha_i \leq 1$, and $x_i \ge 0$ for $i = 1, 2, 3$. A partial CAD (only over the cell $0 <\alpha_1 < 1$, for brevity) is given below. 
    
    \label{ex:full_cad} 
    \begin{figure}[H]
        \centering
        {\scriptsize
\begin{forest}
[{}, phantom
[{$0<\alpha_1<1$}
    [{$\alpha_2=0$}
      [{$\alpha_3=1-\alpha_1$}
        [{$0 \leq x_1 < \frac{1}{\alpha_1}$}
            [{$x_2 \geq 0$}
                [{$x_3 \geq \frac{1-\alpha_1 x_1}{\alpha_3}$}]
            ]
        ]
        [{$x_1 \geq \frac{1}{\alpha_1}$}
            [{$x_2 \geq 0$}
                [{$x_3 \geq 0$}]
            ]
        ]
      ]
    ]
    [{$0<\alpha_2<1-\alpha_1$}
      [{$\alpha_3=1-\alpha_2-\alpha_1$}
        [{$0 \leq x_1 < \frac{1}{\alpha_1}$}
            [{$0 \leq x_2 < \frac{1-\alpha_1 x_1}{\alpha_2}$}
                [{$x_3 \geq \frac{1-\alpha_1 x_1 - \alpha_2 x_2}{\alpha_3}$}]
            ]
            [{$x_2 \geq \frac{1-\alpha_1 x_1}{\alpha_2}$}
                [{$x_3 \geq 0$}]
            ]
        ]
        [{$x_1 = \frac{1}{\alpha_1}$}
            [{$x_2 \geq \frac{1-\alpha_1 x_1}{\alpha_2}$}
                [{$x_3 \geq 0$}]
            ]
        ]
        [{$x_1 > \frac{1}{\alpha_1}$}
            [{$x_2 \geq 0$}
                [{$x_3 \geq 0$}]
            ]
        ]
      ]
    ]
    [{$\alpha_2=1-\alpha_1$}
      [{$\alpha_3=0$}
        [{$0 \leq x_1 < \frac{1}{\alpha_1}$}
            [{$x_2 \geq \frac{1-\alpha_1 x_1}{\alpha_2}$}
                [{$x_3 \geq 0$}]
            ]
        ]
        [{$x_1 > \frac{1}{\alpha_1}$}
            [{$x_2 \geq 0$}
                [{$x_3 \geq 0$}]
            ]
        ]
      ]
    ]
]
]
\end{forest}
}
    \end{figure}
\end{example}

\subsubsection{Extensions: Increasing failure rate}

The \textit{increasing failure rate} (IFR) property is an important property of transition probability matrices in many healthcare and engineering applications \citep{barlow1996mathematical}. It captures the notion that in worse health states, a patient (or machine) is more likely to degrade. More precisely, the rows of the transition probability matrix are in increasing stochastic order. We consider a Markov chain with $\phi$ states, so that we have $\phi$ unit simplices, each with $\phi$ variables: $\{\alpha_{i,j}, i \in [\phi], j \in [\phi] \,|\, \sum_{j=1}^\phi \alpha_{i,j}=1, 0 \le \alpha_{i,j} \le 1, i \in [\phi], j \in [\phi]\}$. 

\begin{definition}[Increasing failure rate (IFR) \citep{alagoz2007determining}]

    Consider $\phi$ unit simplices as above. This set of simplices is said to be IFR if, for each $h \in [\phi]$, we have that $\sum_{\ell = h}^{\phi} \alpha_{i, \ell}$, as a function of $i$, is nondecreasing in $i$, where $i \in [\phi]$.
\end{definition}

The IFR condition introduces a series of linear inequalities as constraints on the simplices. These inequalities span across multiple simplices, so the constraints are not disjoint in the simplex variables anymore. Thus, we cannot use the prior argument to conjunct them to construct the CAD for the set of simplices. 
However, we show that we can still efficiently construct a CAD for a set of simplices with the IFR property. Our strategy will be to construct the full CAD in the usual variable ordering -- $\alpha_{1,1}, \ldots, \alpha_{1,\phi}, \alpha_{2,1}, \ldots, \alpha_{2,\phi}, \ldots, \alpha_{\phi,1}, \ldots \alpha_{\phi, \phi}$. By manipulating the inequalities arising from the IFR condition, we can characterize each $\alpha$ variable in terms of the variables that come before it in this ordering.

\begin{theorem}
    \label{thm:ifr_simplex_cad}
    Consider a set of $\phi$ unit simplices, each with $\phi$ states, $\{\alpha_{i,j}, i \in [\phi], j \in [\phi] \,|\, \sum_{j=1}^\phi \alpha_{i,j}=1, 0 \le \alpha_{i,j} \le 1, i \in [\phi], j \in [\phi]\}$. Also, fix the variable ordering $\alpha_{1,1}, \alpha_{1,2}, \alpha_{1,3}, \cdots, \alpha_{1,\phi}, \alpha_{2,1}, \cdots, \alpha_{2,\phi}, \cdots \alpha_{\phi, \phi}$. If this set of simplices is IFR, their CAD is:
    \begin{itemize}
        \item For each $i$, the cells for $\alpha_{i,j}$, for $j < \phi$ are the cells $0, (0, \sum_{\ell = 1}^j \alpha_{i-1, \ell} - \sum_{\ell = 1}^{j-1} \alpha_{i, \ell}), \sum_{\ell = 1}^j \alpha_{i-1, \ell} - \sum_{\ell = 1}^{j-1} \alpha_{i, \ell}$.
        \item The singular cell for $\alpha_{i,\phi}$ for each $i$ is $\{1 - \sum_{\ell = 1}^{\phi-1} \alpha_{i, \ell}\}$.
        \item If $i=1$, then we replace $\sum_{\ell = 1}^j \alpha_{i-1, \ell}$ with $1$.
        \item The tree propagates any equality constraints on preceding variables in the current cell.
    \end{itemize}
\end{theorem}

\begin{remark}
    The intuition in the bounds on $\alpha$ variables implied by the CAD is that if we write out the $\alpha$ variables in a matrix, a given $\alpha_{i,j}$ is bounded from above by the sum of the $\alpha$ variables in the preceding row (row $i-1$) up to the $j$th column, minus the sum of the $\alpha$ variables in the same row (row $i$) up to the $j-1$th column. So, $\alpha_{i,j}$ is larger if the probabilities in the row above it and to the left are larger, in order to fulfill IFR, and is smaller if the probabilities in the same row but precede it are larger, in order to fulfill row-stochasticity. 
\end{remark}

\begin{example}
    Suppose we have two simplices $\alpha_{1,1} + \alpha_{1,2} = 1$ and $\alpha_{2,1} + \alpha_{2,2} = 1$, and they are IFR. Its CAD is:

    \begin{figure}[H]
        \centering
        {\scriptsize
\begin{forest}
[{}, phantom
[{$\alpha_{1,1} = 0$}
    [{$\alpha_{1,2} = 1$}
        [{$\alpha_{2,1} = 0$}
            [{$\alpha_{2,2} = 1$}]
        ]
    ]
]
[{$0 < \alpha_{1,1} < 1$}
    [{$\alpha_{1,2} = 1 - \alpha_{1,1}$}
        [{$0 \leq \alpha_{2,1} \leq \alpha_{1,1}$}
            [{$\alpha_{2,2} = 1 - \alpha_{2,1}$}]
        ]
    ]	
]
[{$\alpha_{1,1} = 1$}
    [{$\alpha_{1,2} = 0$}
        [{$0 \leq \alpha_{2,1} \leq 1$}
            [{$\alpha_{2,2} = 1 - \alpha_{2,1}$}]
        ]
    ]
]
]
\end{forest}
}
    \end{figure}

    We have merged some cells for brevity. The important difference to note here is how the cells for $\alpha_{2,1}$ depend on $\alpha_{1,1}$, in order to satisfy IFR.
\end{example}

Importantly, the following lemma shows that adding the IFR property does not affect simplex-extensibility of system \eqref{eq:general_form}. 

\begin{lemma}
    \label{lem:ifr_simplex_ext}
    If an instance of \eqref{eq:general_form} is simplex-extensible, then the same system with the additional requirement of IFR remains simplex-extensible. 
\end{lemma}

\subsubsection{Extensions: Other forms of $f^*$}
\label{sec:extensions_fstar}

Thus far, we have focused on $f^*$ that are linear in the $x$-type variables. We now generalize our previous result by allowing for functions of the $x$-type variables, i.e., $f^* = g_0 + \sum_{i=1}^\eta f_i g_i$, where $f_i$ is a polynomial of $x_i$. The main challenge is in keeping track of the projection factors. However, we can generalize it for some polynomials. If, for each $i$, the polynomial $f_i$ has no non-negative roots, i.e., it is sign-invariant over $x_i > 0$, then all the useful properties about simplex-extensibility and projection-definability carry over.

\begin{corollary}
    \label{cor:simplex_ext_f_nonneg}
    Let $f^* = g_0 + \sum_{i=1}^\eta f_i g_i$, where each function $\{f_i\}_{i=1}^\eta$ is a univariate polynomial of $x_i$, and $\{g_i\}_{i=0}^\eta$ is a (possibly multivariate) polynomial of the $\alpha$-type variables. If, for each $i$, $f_i$ has no non-negative roots, then $f^*$ is simplex-extensible if and only if each of $\{g_i\}_{i=0}^\eta$ are sign-invariant over the simplex CAD. Furthermore, if \eqref{eq:general_form} is simplex-extensible with such an $f^*$, it is also projection-definable.
\end{corollary}

Some special cases captured by the preceding corollary are when $f_i$ is a monomial with arbitrary exponent, when $f_i$ has all positive coefficients, or when $f_i$ has all negative coefficients, which is implied by Descartes' rule of signs. There are other such $f_i$ as well, which can be easily identified by root-counting algorithms like Sturm sequences \citep{basu}.\endnote{Writing out the solution formula for such polynomials may introduce difficulties because of having to represent the roots of high-degree univariate polynomials in the CAD tree. This is an issue that should be handled by the computer algebra system so we do not discuss this further here.}



Since polynomials with higher order terms require significant computation to compute projection factors, this corollary is valuable because it allows us to easily check simplex-extensibility for systems where $f^*$ has higher order terms in the $x$-type variables. 

\begin{example}
    Let $f^* = \alpha_1 (x_1^2+x_1^3) + \alpha_2 x_2 + \alpha_3 x_3$, with $\alpha_1+\alpha_2 + \alpha_3=1$. If we tried to check for simplex-extensibility by fully computing the projection factors, we would get $13$ projection factors when eliminating $x_1$, $85$ when eliminating $x_2$, and $751$ when eliminating $x_3$. However, with Lemma \ref{cor:simplex_ext_f_nonneg} we can easily see that $f^*$ is simplex-extensible.
\end{example}

The challenge with generalizing to general univariate polynomials (i.e., with non-negative roots) of $x$-type variables is that, when projecting on a given $x$-type variable, we will generate multiple projection factors in the other $x$-type variables, which makes keeping track of projection factors in subsequent steps onerous. This is unlike the case in Theorem \ref{thm:simplex_ext_xg}, where only a single projection factor in the other $x$-type variables is generated, which keeps the subsequent projection factors simple to compute. Indeed, the technical lemmas in Appendix \ref{sec:some_technical_lemmas} exploit this property. As well, in more general cases, we are not guaranteed projection-definability either, due to the possibility of multiple roots.

The situation is even more difficult if attempting to generalize to multivariate polynomials of the $x$-type variables. This form does not provide any special structure, as indeed any $f^*$ can be written as $f_1g_1 + f_2g_2 + \cdots$, where the $f$ functions are multivariate polynomials of the $x$-type variables and the $g$ functions are multivariate polynomials in the $\alpha$-type variables, simply by taking the individual monomials of $f^*$.

Fortunately, Theorem 2 provides significant generality already, as it applies to all the relevant sensitivity analyses that we introduced in this paper. We now return to that context to demonstrate how our results for system \eqref{eq:general_form} apply to our motivating problem. 


\section{Application to Markov Reward Processes}

As we discussed at the start of Section \ref{sec:specialclass}, system \eqref{eq:general_form} is a general class that contains the polynomial systems associated with our sensitivity analyses of interest. We formalize this result next.


\begin{theorem}
    \label{thm:markov_is_special_system}
    Let $R_\infty$ be defined as in \eqref{eq:infty_reward_nice} and $T \in \mathbb{R}$. Then the system $R_\infty \ge T$ with row-stochastic constraints on $\mat{P}$ and $\mat{\pi}$, and non-negativity constraints on $\mat{r}$:
    
    \begin{enumerate}
        \item is an instance of system \eqref{eq:general_form}
        \item is simplex-extensible, 
        \item is projection-definable, and
        \item has a solution formula that can be constructed via Algorithm \ref{alg:cadsoln_xg}.
    \end{enumerate}
\end{theorem}
\proof{Proof.}


    Letting $f^* := \mat{\pi}^\top \adj{(\mat{I} - \lambda \mat{P})} \mat{r} - T \det{(\mat{I} - \lambda \mat{P})}$, $R_\infty \ge T$ is equivalent to $f^* \ge 0$, so the system is an instance of \eqref{eq:general_form}. The function $f^*$ can be written in the form $g_0 + \sum_{i=1}^\eta x_i g_i$ by setting $g_0 = -T \det{(\mat{I} - \lambda \mat{P})}$, $x_i = r_i$, and $g_i$ as the $i$th entry of $\mat{\pi}^\top \adj{(\mat{I} - \lambda \mat{P})}$. Since $\det{(\mat{I} - \lambda \mat{P})} > 0$ (Lemma \ref{lem:infy_reward_ratio}) and $T$ is fixed, $g_0$ is sign-invariant over the simplex CAD. Similarly, since $\adj(\mat{I} - \lambda \mat{P})$ is component-wise non-negative (Lemma 
    \ref{lem:infy_reward_ratio}), each $g_i$ is non-negative and therefore sign-invariant over the simplex CAD. Thus, the system is simplex-extensible by Theorem \ref{thm:simplex_ext_xg} and projection-definable by Lemma \ref{lem:xg_projdef}, and Algorithm \ref{alg:cadsoln_xg} can be used to construct its CAD solution formula.\halmos
    
    
\endproof

Several remarks follow about the generalizability of this result.

\begin{remark}[Finite-horizon reward]
    The same result holds for the finite-horizon reward with inequality $R_t \ge T$. 
    We set $f^* = R_t - T$, which is of the form $f^* = g_0 + \sum_{i=1}^\eta x_i g_i$ with $\{g_i\}_{i=1}^\eta$ being the entries in $\sum_{m=0}^t \mat{\pi}^\top \lambda^m \mat{P}^m$, per equation \eqref{eq:finite_reward}, which are clearly polynomials in $\mat{\pi}$ and $\mat{P}$. Each of these polynomials are formed by the addition of monomials in $\mat{\pi}$ and $\mat{P}$, and so are non-negative. As well, $g_0 = -T$, which is a constant so is sign-invariant. The results follow by invoking Theorem \ref{thm:simplex_ext_xg} and Lemma \ref{lem:xg_projdef}. 
\end{remark}

\begin{remark}[Other cost-effectiveness quantities]
    By suitably defining $f^*$, Theorem \ref{thm:markov_is_special_system} extends to other sensitivity analyses of interest. 
    
    \begin{hitemize}
        \item For comparing the infinite-horizon rewards of two interventions, labeled $a$ and $b$, the form of the inequality in \eqref{eq:reward_comparison_full} implies that $f^* = \mat{\pi}_a^\top \adj(\mat{I} - \lambda \mat{P}_a) \mat{r}_a \det(\mat{I} - \lambda \mat{P}_b) - \mat{\pi}_b^\top \adj(\mat{I} - \lambda \mat{P}_b) \mat{r}_b \det (\mat{I} - \lambda \mat{P}_a)$, and the $\{g_i\}_{i=1}^\eta$ are the entries in $\mat{\pi}_a^\top \adj(\mat{I} - \lambda \mat{P}_a) \det(\mat{I} - \lambda \mat{P}_b)$ and $- \mat{\pi}_b^\top \adj(\mat{I} - \lambda \mat{P}_b) \det (\mat{I} - \lambda \mat{P}_a)$, which are sign-invariant by similar arguments in the proof above. 
        \item When bounding the NMB, by rearranging equation \eqref{eq:nmb} we set $f^* = \mat{\pi}^\top \adj(\mat{I} - \lambda \mat{P}) (W\mat{b}-\mat{c}) - T \det(\mat{I} - \lambda \mat{P})$. 
    Note that the $x$-type variables are now the elements of $\mat{b}$ and $\mat{c}$.
        \item When bounding the ICER of two interventions $a$ and $b$, we can similarly rearrange \eqref{eq:icer} to obtain a suitable $f^*$, with the $x$-type variables being the costs and benefits of each intervention. 
    \end{hitemize}    
\end{remark}

Lastly, we note that a benefit of taking the CAD perspective is that we can characterize the exact shape of the boundary defined by any of the cost-effectiveness inequalities of interest as a function of $\mat{\pi}$ $\mat{P}$, and/or $\mat{r}$. For example, in \ref{sec:twowaygeom}, we fully characterize the geometry of two-way sensitivity analyses.

\section{Software Implementation}

We developed a Python package, \texttt{markovag}, that allows practitioners to use CAD for sensitivity analyses. The package uses SymPy \citep{sympy}, which is a computer algebra system in Python, to perform symbolic algebraic manipulation. There are two modules.

The first module, \texttt{markovag.markov}, can perform computations with symbolic Markov chains, i.e., where some of the parameters are variables. It contains functions to parse symbolic matrices from input files. One can then perform matrix computations on them. The module comes with functions to symbolically calculate the finite or infinite horizon discounted rewards, as well as the usual health economic metrics of interest like NMB and ICER. In other words, these functions generate the polynomials in equations \eqref{eq:infty_reward_nice} -- \eqref{eq:icer} and their finite horizon analogues. It also provides plotting functions to visualize two- and three-way sensitivity analyses.

The second module, \texttt{markovag.cad},  constructs the CAD. It takes polynomials generated by the first module, appends the necessary stochastic constraints, and builds the CAD tree in accordance with Algorithm \ref{alg:cadsoln_xg}. This development is significant from the larger perspective of computer algebra software, as SymPy currently does not have a CAD solver. In fact, there are only two available implementations of CAD: Mathematica \citep{mathematica_cad}, which is proprietary, and QEPCAD \citep{qepcad}, which is open-source. However, both programs require specific syntax and do not come with specialized methods for cost-effectiveness analysis. Since \texttt{markovag} is built on Python, it is freely available and leverages syntax that is widely used. 

We also note that the two state-of-the-art satisfiability modulo theories solvers, Z3 \citep{z3paper, z3cad} and CVC5 \citep{cvc5paper, cvc5cad} both implement CAD. However, they are satisfiability solvers and so only provide a satisfying point, or indicate that none exist. They do not implement the solution formula construction step. In other words, they only implement Algorithm \ref{alg:caddecision} but not Algorithm \ref{alg:cadsoln}. 

Although \texttt{markovag} only works for instances of system \eqref{eq:general_form} specifically resulting from Markov reward processes, our implementation has many of the fundamentals for constructing a general CAD, including sample point computation and building the tree data structure. Building this out into a full CAD solver is left to future work.

\section{Case Studies}

In this section, we present two case studies to illustrate the use of CAD in sensitivity analysis of Markov reward processes. All results were generated using the \texttt{markovag} package.

\subsection{Synthetic case study}
\label{sec:synth_casestudy}

Consider a three-state Markov chain, where the states represent (1) ``healthy", (2) ``sick", and (3) ``dead". We set the reward of ``healthy" to $r_1$, the reward of ``sick" to $r_2$ and the reward of    ``dead" to $0$. We consider death to be an absorbing state and we compute the infinite horizon expected total reward per equation \eqref{eq:deathstate}. The total reward is a function of the following five parameters: $p_{11}, p_{12}, p_{21}, p_{22}, r_1, r_2$. Observe that they are not a function of $p_{13}$ nor $p_{23}$, i.e., the mortalities of each state, which is implied by equation \eqref{eq:deathstate}.

\subsubsection{Two-way analysis.}

\begin{figure}
    \centering
    \includegraphics[scale=0.6]{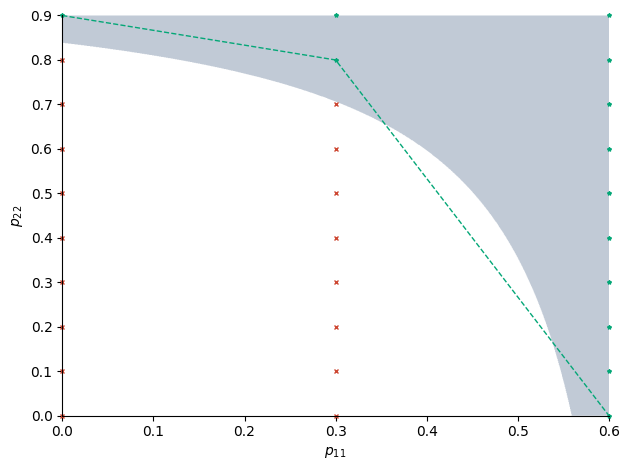}
    \caption{Valid parameter space of a $n=3$ state Markov chain with states (1) healthy, (2) sick, (3) dead, as discussed in Section \ref{sec:synth_casestudy}, with the parameters $p_{12} = 0.4, p_{21} = 0.1, r_1=1, r_2=0.5$. We assert that $R_\infty \geq 3$. Green and red points correspond to the grid search, where green indicates a valid point and red an invalid one. The green lines form the convex hull of the green points.}
    \label{fig:synth}
\end{figure}

First, we conduct a two-way sensitivity analysis of the infinite horizon reward. We fix $p_{12} = 0.4, p_{21} = 0.1, r_1=1, r_2=0.5$, and let $p_{11}$ and $p_{22}$ vary. Given the fixed parameters, we have the implied bounds $p_{11} \in [0, 0.6]$ and $p_{2,2} \in [0, 0.9]$. We wish to find the values of $p_{11}$ and $p_{22}$ that guarantee $R_\infty \geq 3$. Using \texttt{markovag.markov}, we symbolically compute the expected total reward, form the inequality bounding it, and plot the resulting valid parameter space (Figure \ref{fig:synth}). 

Notice that this parameter space is \textit{not convex}. Without CAD, the traditional approach of using a fixed parameter grid would require shrinking the grid size to get a good approximation of the boundary. However, it would not be known a priori how small the grid needs to be to achieve a desired error. Our exact approach circumvents this issue.

In Figure \ref{fig:synth}, we overlay a grid of lattice points, where green points satisfy the inequality and red points do not. The naive approach of identifying neighboring grid points that satisfy the inequality (i.e., cost-effective) may lead to an incorrect conclusion that all convex combinations of those parameter values also satisfy the inequality. 

Another result is that at $p_{11} = 0.3$, we can tell from the grid that $p_{22} = 0.8$ the inequality is satisfied by at $p_{22} = 0.7$ it is not. So the traditional analysis only knows that the cost-effectiveness boundary is crossed somewhere in the range $[0.7,0.8]$, whereas the CAD analysis clearly shows that the boundary is very close to $0.7$.

\subsubsection{Multi-way analysis via CAD}

To illustrate the power of the CAD approach, we consider an eight-way sensitivity analysis. We allow all of $p_{11}, p_{12}, p_{13}, p_{21}, p_{22}, p_{23}, r_1, r_2$ to vary. The difficulty in searching and visualizing this high-dimensional grid means that an eight-way analysis is most likely never done in practice. 
We once again assert that $R_\infty \geq 3$. Using the \texttt{markovag} package we construct the CAD tree and restrict our attention to the full-dimensional cells (cells which are intervals) for the probabilities. We explore the effect of the probabilities on the permissible rewards, as well as the interaction between the two rewards $r_1$ and $r_2$. For example, we find that if 
\begin{equation*}
    r_1 \geq \frac{3 p_{1 1} p_{2 2} - 3 p_{1 1} - 3 p_{1 2} p_{2 1} - 3 p_{2 2} + 3}{1 - p_{2 2}}
\end{equation*}
then $r_2 \ge 0$.  This result is useful, for example, if the reward for being in the healthy state is well-known and accepted but the reward for being in the sick state is highly uncertain and patient dependent. A policymaker only needs to know that if a patient receives a sufficiently high reward for being in the healthy state, then any reward value in the sick state would still lead to a cost-effective result. Furthermore, if the reward in state 1 can be written as benefit minus cost, a lower bound on $r_1$ implies an upper bound on cost, such that if the policymaker can bring the cost of a hypothetical therapeutic below the bound, then the intervention would be considered cost-effective even without knowing the exact reward for being in state 2.

Another benefit of our approach is that it immediately elucidates the impact of adding or removing restrictions on the parameter values. For example, the analysis up to now did not include the IFR condition. 
Without this condition, the CAD analysis shows that, e.g., $p_{11} = 0.2$, $p_{12} = 0.5$, $p_{21} = 0.3$, $p_{22} = 0.5$, $r_1 \geq 1.5$, $r_2 \geq 0$ results in $R_\infty \ge 3$, suggesting that this combination of parameters is associated with a cost-effective intervention. However, with the IFR condition imposed, the CAD tree would no longer contain this set of parameter values due to the additional restrictions imposed on $p_{21}$ and $p_{22}$. Lowering $p_{21}$ to satisfy IFR, e.g., $p_{21} = 0.1$, would result in $r_1 \geq 2.15$, $r_2 \geq 0$. This result implies that a higher reward is needed in the healthy state to retain cost-effectiveness. 

\subsection{Drones for cardiac arrest response} 

In this subsection, we use our CAD approach to re-analyze an existing cost-effectiveness analysis from the literature \citep{dronecea}. This paper studied the cost-effectiveness of using drones to deliver automated external defibrillators (AEDs) to the patients suffering out-of-hospital cardiac arrest (OHCA). Rapid response is critical for survival from OHCA and drones as an AED delivery vector have received significant attention in recent years, including network modeling \citep{boutilier2017optimizing,boutilier2022drone}, to dispatching \citep{chu2021machine}, to pilot tests \citep{cheskes2020improving, claesson2017time} to real-world implementation \citep{schierbeck2022use, schierbeck2023drone}. However, despite rapid progress and growing interest in implementation, cost-effectiveness of this intervention has received limited attention. 

The analysis in \citet{dronecea} used data on 22,017 real OHCAs in Ontario, Canada to evaluate the potential cost-effectiveness of 1000 different drone network configurations (i.e., the interventions) designed using various optimization models. A decision tree compared drone-augmented ambulance response against the baseline of ambulance-only response. If the drone arrives before the ambulance, a responder would apply the AED with some probability. If the AED is applied, machine learning models were used to predict whether the patient survived or not, and if survived, their neurological status as measured by the modified Rankin scale (mRS), where 0 is perfect health and 6 is death. If the AED is not applied or the drone arrives after the ambulance, then the patient follows the trajectory recorded in the historical data. Finally, a Markov model was used to simulate post-arrest trajectories and calculate accumulated QALYs and costs over the remaining patient lifetime. Transitions between health states, costs and utilities were mRS-specific and obtained from the medical literature. 
We refer the reader to \citet{dronecea} for full details on the methods. As drones remain a novel intervention with limited real-world implementation, several important parameters for the analysis were not well-established in the literature, and so numerous one-way sensitivity analyses were performed. 

Now, we perform a much more thorough sensitivity analysis using the tools developed in this paper and our software. We
re-implemented the decision tree and Markov model from \citet{dronecea} using \texttt{markovag.markov}. We consider one particular drone network as an example, a 20-drone network that maximized 5-minute coverage of cardiac arrests, as this was the smallest drone network found to be cost-effective. We performed a three-way sensitivity analysis on the following three key parameters:
\begin{itemize}
    \item $m$: the drone cost multiplier. \citet{dronecea} obtained estimates of the cost to operate the drone network from two drone companies. However, as the technology evolves and the complexity of operating a large network changes (currently, the largest operating network in the world, in Sweden, has only five drones), this cost may also change in unexpected ways. Thus, we multiplied the estimated drone costs by different values of $m$. In the original paper, values of $m \in \{2, 5\}$ were tested.
    \item $p_u$: the probability of the drone-delivered AED being used. A value of 0.457 was used for this probability, obtained from a prior study on bystander use of static AEDs. However, there remains significant uncertainty about the true usage probability and how it might vary across geographies. Thus, in the original one-way sensitivity analysis, a range of $[0.05, 0.75]$ was tested. 
    \item $p_0$: the first-year mortality of a patient discharged with mRS 0. In the original paper, this parameter was found to be one of the most influential on NMB. It had a default value of $0.013$, and in the original one-way sensitivity analysis, a range of $[0.002, 0.048]$ was tested. 
\end{itemize}

We used \texttt{markovag.markov} to symbolically calculate the difference in NMB between the drone case and the no-drones case. Then we used \texttt{markovag.cad} to construct the CAD for the inequality representing ``drone NMB $\ge$ status quo NMB''. The most important cell in the CAD is
\begin{align*}
    0 < p_{u} < 1 
    \quad \land  \quad 0 < p_0 < 1 
    \quad \land \quad 0 < m \leq -132.46 p_{0} p_{u} + 154.46 p_{u} 
\end{align*}
We omit the other cells as they represent extreme cases that are uninteresting or unrealistic for the policymaker (e.g., $p_u = 1$, meaning universal drone AED use). Some insights can be gleaned by analyzing the first derivatives of the multilinear function bounding $m$: 
\begin{itemize}
    \item With respect to $p_u$: we have $ -132.46 p_{0} + 154.46$, which is strictly positive for $0 < p_0 < 1$, so the bound on $m$ is strictly increasing in $p_u$. As expected, with a higher probability of drone AED use, higher drone costs are acceptable due to the larger benefits.
    \item With respect to $p_0$: we have $-132.46 p_{u}$ which is strictly negative for $0 < p_u < 1$, so the bound on $m$ is strictly decreasing in $p_0$. This result is not easy to obtain otherwise and is somewhat unexpected. A higher probability of death leads to lower utilities obviously, but also lower costs, as there are non-trivial costs associated with treating surviving patients. In this case, the analysis reveals that the increased costs outweigh the increased benefits, so that with a higher probability of death, the acceptable level of drone costs \textit{decreases}.
\end{itemize}

From this CAD representation, a policymaker can trace a path down the tree to determine the validity of a set of parameters.

For illustrative purposes, we visualize the parameter regime over which the drone network is cost-effective, found by \texttt{markovag}, by fixing $m \in \{1,5\}$ and letting $p_0$ and $p_u$ vary. For comparison, we consider the usual one-way and two-way sensitivity analyses using a grid with five evenly spaced points, $\{0, 0.25, 0.5, 0.75, 1\}$. For the one-way analysis, we fixed one parameter at the default value according to \citet{dronecea} ($p_0 = 0.013, p_u = 0.457$) and tested the free parameter at the five points above. For the two-way sensitivity analysis, we construct a 5x5 grid of the two parameters. Grid points where the inequality holds (does not hold) are denoted by a green star (red cross). 

\begin{figure}
    \centering
    \includegraphics[scale=0.6]{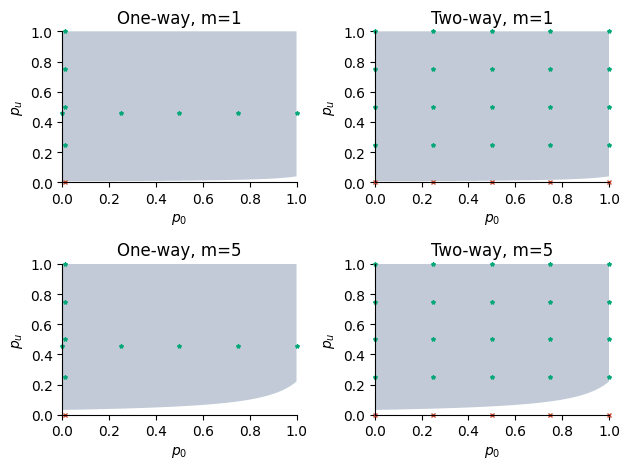}
    \caption{Visualization of the parameter space over which the drone network is cost-effective. The shaded gray region is the exact analytic solution obtained by \texttt{markovag}. The points represent the traditional mesh grid approach, with green representing a valid point and red invalid.}
    \label{fig:mainplot}
\end{figure}

The plots in Figure \ref{fig:mainplot} show that the range of parameter values associated with cost-effectiveness is quite large, much more than suggested by the one-way or even two-way analyses, and that the nonlinear boundary is easily identified. Methodologically, a grid search will always miss some part of the true cost-effective region and it remains difficult to know a priori how fine a grid is needed to approximate well the boundary. In the context of the drone application, these results elucidate that the probability of drone AED usage may be much lower than expected while still leading to a cost-effective intervention. Also, sensitivity of the cost-effectiveness finding to $p_0$ seems to be more apparent only when considering the interaction with other parameters. 


\section{Conclusion}

In this paper, we studied the problem of performing sensitivity analysis on metrics that are derived from Markov reward processes. We were motivated specifically by cost-effectiveness analyses in healthcare. We showed that every such analysis has an equivalent semialgebraic representation. Our framework can encompass other common extensions in the healthcare literature, like the IFR assumption. We then showed that these systems have a special structure which allows us to construct their CADs in a simpler way. We call this property simplex-extensible, and it is more general than just the polynomials arising from our study of Markov reward processes. Lastly, we developed software which allows practitioners to use our approach, and we demonstrated that it finds regions of cost-effectiveness that are missed by the classical grid search method.

\theendnotes
	
\newpage 

\bibliographystyle{informs2014} 
\bibliography{ref}

\ACKNOWLEDGMENT{This paper is supported by the Natural Sciences and Engineering Research Council of Canada (NSERC).}

\newpage

\ECSwitch

\ECHead{Electronic Companion}

\section{Background on Hong Projection Operator}
\label{sec:hong}

Here, we describe Hong's projection operator \citep{hong1990improvement}. First, we introduce some new notation. To use a projection operator, we have to choose a main variable, or mvar. The polynomials will be treated as univariate polynomials in the mvar. Then, for a polynomial $f$, the leading term, $\ldt(f)$, is the term with the highest exponent of the mvar, the leading coefficient, $\ldcf(f)$, is the coefficient of $\ldt(f)$, and the degree, $\deg(f)$ is the highest exponent with which the mvar appears.

Two crucial concepts in algebraic geometry are the \textit{reducta} and the \textit{principal subresultant coefficients}, which we define below.

\begin{definition}[Reducta]
    The reductum of a polynomial $f$ with a chosen mvar is $\red(f) = f-\ldt(f)$. We inductively define the $i$th level reductum as: $\red^0(f) = f$ and $\red^i(f) = \red(\red^{i-1}(f))$. Then, the reducta set $\redset(f) = \{\red^i(f), 0 \leq i \leq \deg(f), \red^i(f) \neq 0\}$.
\end{definition}

\begin{example}
    Let $f=x_1 x_2 + x_3 x_4 -1$, with mvar $x_1$. Then, $\red^0 = x_1 x_2 + x_3 x_4 -1$, and $\red^1 = x_3 x_4 - 1$. So, $\redset(f) = \{x_1 x_2 + x_3 x_4 -1, x_3 x_4 - 1\}$.
\end{example}

\begin{definition}[Principal subresultant coefficients]

Let there be two polynomials $f,g$, with degrees $p,q$, respectively, having chosen the mvar $x$. Let $p>q$ (resp. $p=q$), and fix $0 \leq i \leq q$ (resp. $0 \leq i \leq p-1$). Define the \textit{$i$th Sylvester-Habicht matrix}, denoted $\sylvha_i(f,g)$, as the matrix whose rows are $x^{q-i-1}f, x^{q-i-2}f, \cdots, f, g, \cdots, x^{p-i-1}g$, considered as vectors in the basis $[x^{p+q-i-1}, \cdots, x, 1]$; it has $p+q-i$ columns and $p+q-2i$ rows. Then, the $i$th principal subresultant coefficient, denoted $\psc_i$, is the determinant of the submatrix of $\sylvha_i(f,g)$ obtained by taking the first $p+q-2i$ columns. Then, the PSC set $\pscset(f,g) = \{\psc_i(f,g), 0 \leq i \leq \min(\deg(f), \deg(g)), \psc_i(f,g) \neq 0\}$.
\end{definition}

\begin{example}
    Let $f =3x^2 + 5x +6, g=4x^2+2x+1$, with mvar $x$, and say we want $\psc_0$. The Sylvester-Habicht matrix will have $4$ columns and $4$ rows. The basis for the rows is $[x^3, x^2, x, 1]$. For the first row, we compute $x \cdot f = 3x^3 + 5x^2 + 6x$, which is $[3,5,6,0]$ in our basis. The second row is $[0,3,5,6]$. For the third row, we take $g$, which yields $[0,4,2,1]$. For the fourth row, we compute $x \cdot g$, yielding $[4,2,1,0]$ in our basis. Hence, we have the following matrix: 
    \begin{equation*}
    \begin{bmatrix}
    3 & 5 & 6 & 0 \\
    0 & 3 & 5 & 6 \\
    0 & 4 & 2 & 1 \\
    4 & 2 & 1 & 0
    \end{bmatrix}
    \end{equation*}
    
    Then, $\psc_0$ is the determinant of the whole matrix, which is $-343$.
\end{example}

The operations defined above are used in the Hong projection operator. We let $\mathcal{D}$ represent the derivative operator. Then, the Hong projection operator $\projh(F)$ of a set $F$ of polynomials is \citep{hong1990improvement}:

\begin{align}
    &\projh(F) = \projone(F) \cup \projtwo(F) \\
    &\projone(F) = \bigcup_{\substack{f \in F \\ f' \in \redset(f)}} \left[ \{\ldcf(f')\} \cup \pscset(f', \mathcal{D}f')\right] \\
    &\projtwo(F) = \bigcup_{\substack{f,g \in F \\ f \prec g}} \bigcup_{\substack{f' \in \redset(f)}} \pscset(f',g)
\end{align}

Above, $f \prec g$ denotes an arbitrary linear ordering, to not loop over redundant pairs. The operator $\projone$ tells us to loop over the set $F$, and loop over the reducta set for each, and apply $\ldcf$ and $\pscset$. The operator $\projtwo$ tells us to loop over (non-redundant) pairs $f,g \in F$, take the first one's reducta set, and calculate the PSC set between each of them with $g$. \citet{hong1990improvement} proved that this is a valid projection operator.

Now that we have defined the Hong projection operator, we demonstrate an example, using the sphere, as in Examples \ref{ex:sphere_cad} and \ref{ex:sphere_cad_execution}. Recall in Example \ref{ex:sphere_cad_execution} we skipped over the actual execution of the projection phase, but we give it here below.

\begin{example}
    Suppose we have the polynomial $x^2+y^2+z^2-1$. We want to compute the projection factors of this system, which only has a single polynomial. In keeping with our notation, we set $F_3 = \{x^2+y^2+z^2-1\}$. We now eliminate $z$ and then $y$.

    \begin{henumerate}
        \item Eliminate $z$.
        \begin{hitemize}
            \item Apply $\projone$.
                \begin{hitemize}
                    \item Set $f = x^2+y^2+z^2-1$. Then, $\redset(f) = \{x^2+y^2+z^2-1, x^2+y^2-1\}$.
                    \begin{hitemize}
                        \item Set $f' = x^2+y^2+z^2-1$. We have $\ldcf(f') = 1$. Also, $\mathcal{D} f' = 2z$. We now compute $\pscset$. 
                        \begin{hitemize}
                            \item Set $i=0$. $\sylvha_0$ is a $3 \times 3$ matrix, who's rows have the basis $[z^2, z, 1]$. We have:
                            \begin{equation*}
                                \sylvha_0(f', \mathcal{D} f') = 
                                \begin{bmatrix}
                                    1 & 0 & x^2 + y^2 - 1 \\
                                    0 & 2 & 0 \\
                                    2 & 0 & 0
                                \end{bmatrix}
                            \end{equation*}
                            Then, $\psc_0(f', \mathcal{D} f')$ is the determinant of the entire above $\sylvha$ matrix, equalling $4 - 4x^2 - 4y^2$.

                            \item Set $i=1$. $\sylvha_1$ is a $1 \times 2$ matrix, who's rows have the basis $[z, 1]$. We have:
                            \begin{equation*}
                                \sylvha_1(f', \mathcal{D} f') = 
                                \begin{bmatrix}
                                    2 & 0
                                \end{bmatrix}
                            \end{equation*}
                            Then, $\psc_1(f', \mathcal{D} f')$ is the determinant of the submatrix of the above matrix taking the first $1$ columns, equalling $2$.
                        \end{hitemize}

                        \item Set $f' = x^2+y^2-1$. We have $\ldcf(f')=x^2+y^2-1$. Also, $\mathcal{D} f' = 0$. We now compute $\pscset$. As both are degree 0 in $z$, it is empty.
                    \end{hitemize}
                \end{hitemize}
            \item Apply $\projtwo$: can skip as we only have a single polynomial.
            \item Hence, $F_2 = \{1, 2, x^2+y^2-1, -4x^2-4y^2+4\}$. We can drop the constants (see Appendix \ref{sec:some_technical_lemmas}). 
        \end{hitemize}

        \item Eliminate $y$.
        \begin{hitemize}
            \item Apply $\projone$.
                \begin{hitemize}
                    \item Set $f = x^2+y^2-1$. Then, $\redset(f) = \{x^2+y^2-1, x^2-1\}$.
                    \begin{hitemize}
                        \item Set $f' = x^2+y^2-1$. We have $\ldcf(f') = 1$. Also, $\mathcal{D} f' = 2y$. We now compute $\pscset$. We omit some details as we showed a detailed run-through above.
                        \begin{hitemize}
                            \item Set $i=0$. We have:
                            \begin{equation*}
                                \sylvha_0(f', \mathcal{D} f') = 
                                \begin{bmatrix}
                                    1 & 0 & x^2 - 1 \\
                                    0 & 2 & 0 \\
                                    2 & 0 & 0
                                \end{bmatrix}
                            \end{equation*}
                            Then, $\psc_0(f', \mathcal{D} f')$ is the determinant of the entire above $\sylvha$ matrix, equalling $4 - 4x^2$.

                            \item Set $i=1$. We have:
                            \begin{equation*}
                                \sylvha_1(f', \mathcal{D} f') = 
                                \begin{bmatrix}
                                    2 & 0
                                \end{bmatrix}
                            \end{equation*}
                            Take the first $1$ columns: $\psc_1(f', \mathcal{D} f')$ equals $2$.
                        \end{hitemize}

                        \item Set $f' = x^2-1$. We have $\ldcf(f')=x^2-1$. Also, $\mathcal{D} f' = 0$. Lastly, $\pscset(f', \mathcal{D} f') = \emptyset$.
                    \end{hitemize}

                    \item Set $f = -4x^2-4y^2+4$. Then, $\redset(f) = \{-4x^2-4y^2+4, -4x^2+4\}$.
                    \begin{hitemize}
                        \item Set $f' = -4x^2-4y^2+4$. We have $\ldcf(f') = -4$. Also, $\mathcal{D} f' = -8y$. We now compute $\pscset$. 
                        \begin{hitemize}
                            \item Set $i=0$. We have:
                            \begin{equation*}
                                \sylvha_0(f', \mathcal{D} f') = 
                                \begin{bmatrix}
                                    -4 & 0 & 4-4x^2 \\
                                    0 & -8 & 0 \\
                                    -8 & 0 & 0
                                \end{bmatrix}
                            \end{equation*}
                            Then, $\psc_0(f', \mathcal{D} f') = -256+265x^2$.

                            \item Set $i=1$. We have:
                            \begin{equation*}
                                \sylvha_1(f', \mathcal{D} f') = 
                                \begin{bmatrix}
                                    -8 & 0
                                \end{bmatrix}
                            \end{equation*}
                            Take the first $1$ columns: $\psc_1(f', \mathcal{D} f') = -8$.
                        \end{hitemize}

                        \item Set $f' = -4x^2+4$. We have $\ldcf(f')=-4x^2+4$. Also, $\mathcal{D} f' = 0$. Lastly, $\pscset(f', \mathcal{D} f') = \emptyset$.
                    \end{hitemize}
                \end{hitemize}
                
            \item Apply $\projtwo$. 
                \begin{hitemize}
                    \item Set $f=x^2+y^2-1, g=-4x^2-4y^2+4$. We have $\redset(f) = \{x^2+y^2-1, x^2-1\}$.
                    \begin{hitemize}
                        \item Set $f' = x^2+y^2-1$. We have $\psc_0 = 0$ and $\psc_1 = 0$ (we omit the construction of the $\sylvha$ matrices).
                        \item Set $f' = x^2-1$. We have $\psc_0 = 1 - 2x^2 + x^4$.
                    \end{hitemize}
                \end{hitemize}
            
            \item Hence, $F_1 = \{-8, -4, 0, 1, 2, -256+256x^2, -4x^2+4, x^2-1, x^4-2x^2+1\}$. We can drop the constants (see Appendix \ref{sec:some_technical_lemmas}). 
        \end{hitemize}
    \end{henumerate}

    Therefore, we have the projection factors after eliminating $z$: $F_2= \{x^2+y^2-1, -4x^2-4y^2+4\}$ and the projection factors after eliminating $y$: $F_1= \{-256+256x^2, -4x^2+4, x^2-1, x^4-2x^2+1\}$.
\end{example}

\section{Technical Lemmas About Reducta and Principal Subresultant Coefficients}
\label{sec:some_technical_lemmas}

We introduce several technical lemmas about reducta and principal subresultant coefficients through which we will study the behaviour of the Hong projection for instances of system \eqref{eq:general_form}.

First, we have the following lemmas about reducta sets.

\begin{lemma}
    \label{lem:red_dg0}
    If $f$ is degree $0$, then $\bigcup_{f' \in \redset(f)} \left[ \ldcf(f') \cup \pscset(f', \mathcal{D}f') \right] = \{f\}$.
\end{lemma}
\proof{Proof.}
    $\redset(f) = \{\red^0(f)\} = f$, and $\ldcf(f) = f$, and $\mathcal{D}f=0$, so $\pscset(f,\mathcal{D}f)=\emptyset$. Hence, we only have $\{f\}$. \halmos
\endproof

\begin{lemma}
    \label{lem:red_dg1}
    If $f$ is degree $1$, then $\bigcup_{f' \in \redset(f)} \left[ \ldcf(f') \cup \pscset(f', \mathcal{D}f') \right] = \{\ldcf(f), f - \ldt(f)\}$.
\end{lemma}
\proof{Proof.}
    Here, $\redset(f) = \{f, f-\ldt(f)\}$. Now, $\mathcal{D}f = f - \ldt(f)$ is degree-zero, so to compute $\pscset(f,\mathcal{D}f) = \{\psc_0(f, \mathcal{D}f)\}$, note that $\sylvha_0(f, \mathcal{D}f)$ is the $1 \times 1$ matrix with the entry $f-\ldt(f)$, so the determinant and hence $\psc_0(f, \mathcal{D}f) = f-\ldt(f)$. Lastly, we will include $\ldcf(f)$, and $\ldcf(f-\ldt(f))=f-\ldt(f)$, because it is degree-zero. So we are left with $\{\ldcf(f), f-\ldt(f)\}$. \halmos
\endproof

Next, we have the following lemmas about principal subresultant coefficients.

\begin{lemma}
    \label{lem:psc_dg0}
    Let $f,g$, with $g$ degree 0 and $\deg(f)=d>0$. Then $\pscset(f,g)=\{(-1)^{\lfloor d/2 \rfloor} \cdot g^{d}\}$.
\end{lemma}
\proof{Proof.}
    Here, $\pscset(f,g) = \{\psc_0(f,g)\}$. We form $\sylvha_0(f,g)$ as follows: it is a $d \times d$ matrix, and no rows will correspond to $f$, because $\deg(g)=0$. So, it is the antidiagonal matrix with entries $g$. We perform $\lfloor d/2 \rfloor$ row swaps (i.e., first and last, second and second last, etc.) to make it a diagonal matrix, noting that every swap causes a sign change in the determinant. Then, the determinant of the diagonal matrix is $g^{d}$, which is multiplied by $(-1)^{\lfloor d/2 \rfloor}$ because of the swaps. \halmos
\endproof

\begin{lemma}
    \label{lem:psc_dg0_next}
    Let $f,g$, with $g$ degree 0, and $f$ have a degree sequence, i.e., exponents of $x$ with a non-zero coefficient, of $D$. Then, $\bigcup_{\substack{f' \in \redset(f)}}\pscset(f',g)= \{(-1)^{\lfloor d/2 \rfloor} \cdot g^{d}, d \in D \}$. 
\end{lemma}
\proof{Proof.}
    The degrees of the reducta set of $f$ are exactly the degree sequence $D$. The result follows from Lemma \ref{lem:psc_dg0}. \halmos
\endproof

\begin{remark}
    For the purposes of CAD, the sign change, though mathematically accurate, is not necessary to keep track of, because ultimately we care about the roots of the projection factors, which is unaffected by negation. Therefore, when we apply Lemmas \ref{lem:psc_dg0} and \ref{lem:psc_dg0_next}, we will ignore the  $(-1)^{\lfloor d/2 \rfloor}$ sign.
\end{remark}

We summarize some key takeaways here:

\begin{itemize}
    \item For $\projone$, degree $0$ polynomials (in the mvar, although they may still have positive degree in other variables), are just copied over into the projection factor set.
    \item For $\projtwo$, we can ignore pairs with both degree-zero polynomials. If the pair has a single degree-zero polynomial, we can easily use Lemma \ref{lem:psc_dg0_next}.
    \item For constants, i.e., degree 0 in all variables, we need not store them, as they will be propagated to every subsequent projection factor set, and in the base phase do not contribute to creating any cells.
\end{itemize}

\section{Omitted Proofs from Section \ref{sec:cad}}
\label{sec:omitted_proofs}

\proof{Proof of Theorem \ref{thm:simplex_cad}}
    Observe that a valid representation of a simplex can be constructed as follows: $\alpha_1 \in [0,1]$, $\alpha_2 \in [0, 1-\alpha_1]$, $\alpha_3 \in [0, 1-\alpha_1 - \alpha_2]$, and so on, until $\alpha_\tau = 1 - \sum_1^{\tau-1} \alpha_i$. Clearly, each variable is in $[0,1]$. The sum $\sum_1^\tau \alpha_i \geq 1$, by summing the lower bounds of these intervals. Next, by looking at the partial sums of the upper bounds, we can see that $\alpha_1 + \alpha_2 \leq 1$, $\alpha_1 + \alpha_2 + \alpha_3 \leq 1$, and so on, so that $\sum_1^\tau \alpha_i \leq 1$. Therefore, $\sum_1^\tau \alpha_i = 1$, as required by the simplex constraint.

    Then, we convert this into a valid CAD by decomposing each of the defining intervals into their cells, namely the endpoints and the open interval. \halmos
\endproof

\proof{Proof of Corollary \ref{cor:gluing_simplex_cad}}
    We proceed by induction. In the base case, with one simplex, it is trivially true. For the inductive step, assume we have the CAD of the conjunction of $m$ simplices. Then, to lift this CAD to include the $m+1$-st simplex, first note that because the $\alpha$-type variables only appear in a single simplex, the lifting is identical over each cell, and indeed the lifting process will simply yield the additional $m+1$-st simplex itself. Therefore, each cell in the CAD of the conjunction of $m+1$ simplices is the conjunction of a single cell from each of the simplices individually, and the full CAD is the disjunction over all such cells, as needed.
    
    Note that this also follows from the distributivity of conjunction over disjunction for an arbitrarily indexed family of sets: e.g., see \citet[Theorem 5.21]{monk}. \halmos
\endproof

\proof{Proof of Corollary \ref{cor:simplex_cad_numcells}}
    The CAD of a single simplex has at most 3 children at each level, except at the last level where it definitely has one. So the number of cells for the $i$th simplex is $O(3^{\tau_i-1}) = O(3^{\tau_i})$. Then, if we conjunct $\phi$ simplices together, we multiply the number of cells: $O(3^{\tau_1} \times 3^{\tau_2} \cdots \times 3^{\tau_\phi}) = O(3^{\sum_{i=1}^\phi \tau_i})$. Since each $\tau_i \leq \max_{1 \leq j \leq \phi} \tau_j$, we also have the bound $O(3^{\phi \cdot \max_j \tau_j})$. \halmos
\endproof

\proof{Proof of Theorem \ref{thm:ifr_simplex_cad}}
    First, from the definition of IFR, we can write that, for a given $h \in [1, \phi]$, we have:

    \begin{align*}
        & \sum_{\ell = h}^{\phi} \alpha_{1,\ell} \leq \sum_{\ell = h}^{\phi} \alpha_{2,\ell} \leq 
        \cdots \leq 
        \sum_{\ell = h}^{\phi} \alpha_{\phi,\ell} \\
        & \iff 1 - \sum_{\ell = 1}^{h-1} \alpha_{1,\ell} \leq
        1 - \sum_{\ell = 1}^{h-1} \alpha_{2,\ell} \leq 
        \cdots \leq
        1 - \sum_{\ell = 1}^{h-1} \alpha_{\phi,\ell} \\
        & \iff \sum_{\ell = 1}^{h-1} \alpha_{1,\ell} \geq
        \sum_{\ell = 1}^{h-1} \alpha_{1,\ell} \geq 
        \cdots \geq 
        \sum_{\ell = 1}^{h-1} \alpha_{1,\ell}
    \end{align*}

    Note that when $h=1$, the inequalities are just $1 \leq \cdots \leq 1$, which adds no information.
    
    For the $\alpha_{i,j}$ variables for $i=1$, we have the usual characterization that $\alpha_{1,1} \in [0,1]$, $\alpha_{1,2} \in [0, 1 - \alpha_{1,1}]$, and so on, until $\alpha_{1,\phi} = 1 - \alpha_{1, \phi-1} - \alpha_{1, \phi-2} - \cdots - \alpha_{1,1}$.

    For $\alpha_{i,j}$, $i > 1$, from the inequalities in the IFR definition when $h=j+1$, we get that:

    \begin{align*}
        & \sum_{\ell = 1}^j \alpha_{i-1, \ell} \geq \sum_{\ell = 1}^{j} \alpha_{i, \ell} \\
        & \iff \alpha_{i,j} \leq \sum_{\ell = 1}^j \alpha_{i-1, \ell} - \sum_{\ell = 1}^{j-1} \alpha_{i, \ell}
    \end{align*}

    Recall that the usual bound $\alpha_{i,j} \leq 1 - \sum_{\ell = 1}^{j-1} \alpha_{i, \ell}$. Since $\sum_{\ell = 1}^j \alpha_{i-1, \ell} \leq 1$, the above is a stricter bound than we would have without the IFR constraint. Also, this provides an exact characterization of $\alpha_{i,j}$ in terms of the $\alpha$-variables that precede it in our variable ordering.

    Lastly, as implied by row-stochasticity, we must have $\alpha_{i,\phi} = 1 - \sum_{\ell = 1}^{\phi-1} \alpha_{i, \ell}$. The CAD decomposition follows. \halmos
\endproof

\proof{Proof of Theorem \ref{thm:simplex_ext_xg}}
    We eliminate the $x$-type variables in turn. First, if we set $x_\eta$ as the mvar, then for $\projone$, we copy over the simplex constraints, and for $f^*$, we keep the $\ldcf$ (which is $g_\eta$) and subtract the $\ldt$ (which is $x_\eta g_\eta$). For $\projtwo$, we raise the simplex constraints to the power of one, but we already have them. So, this yields a projection factor set of $\{g_\eta, g_0 + \sum_1^{\eta-1} x_i g_i\}$ plus the simplex constraints. 

    We then similarly eliminate $x_{\eta-1}$, which yields the projection factor set $\{g_{\eta}, g_{\eta-1}, g_0 + \sum_1^{\eta-2} x_i g_i\}$ plus the simplex constraints. Inductively, after eliminating all $x$-type variables, we are left with $\{g_i\}_0^\eta$ plus the simplex constraints.

    By definition, the instance is simplex-extensible if and only if each of these projection factors is sign-invariant over the simplex CAD. Trivially, the simplex constraints are sign-invariant over the simplex. Therefore, it is simplex-extensible if and only if each of $\{g_i\}_0^\eta$ are sign-invariant over the simplex CAD. \halmos
\endproof

\proof{Proof of Lemma \ref{lem:ifr_simplex_ext}}
    Each cell of the simplex with the IFR condition is a subset of the simplex without the IFR condition, so that if a polynomial is sign-invariant over a cell in the non-IFR simplex CAD then it is also sign-invariant in the corresponding IFR simplex CAD. \halmos
\endproof

\proof{Proof of Lemma \ref{lem:xg_projdef}}
    First, we already have the solution formula for the simplex CAD, so we need only focus on the solution formula when lifting to the $x$-type variables. The set of all projection factors, following the proof of Theorem \ref{thm:simplex_ext_xg}, is:
    \begin{align*}
        &\left\{g_\eta, g_0 + \sum_1^{\eta-1} x_i g_i\right\} \bigcup \left\{g_\eta, g_{\eta-1}, g_0 + \sum_1^{\eta-2} x_i g_i\right\} \bigcup \cdots \bigcup \left\{g_i \right\}_0^\eta \\
        &=  \left\{g_i \right\}_0^\eta \bigcup \left\{ g_0 + \sum_{j=1}^{\eta-i} x_i g_i \right\}_{i=1}^{\eta-1}
    \end{align*}

    We can construct a solution formula from the projection factors as we discussed previously, by conjuncting atoms in accordance with the signs of these projection factors. As $\{g_i\}_0^\eta$ are sign-invariant over each cell, due to simplex-extensibility, it is redundant to include them. Secondly, atoms formed from the following set:

    \begin{equation*}
    \left\{ g_0 + \sum_{j=1}^{\eta-i} x_i g_i \right\}_{i=1}^{\eta-1}
    \end{equation*}

    are linear in the respective $x_i$, hence have at most one root over the cell that we are lifting from, so that the sign of the projection factor uniquely determines a region. Therefore, this system is projection-definable. \halmos
\endproof

\proof{Proof of Corollary \ref{cor:full_cad_size}}
    The number of cells of the simplex CAD is as per Corollary \ref{cor:simplex_cad_numcells}. After the simplex CAD, each level has at most 4 cells, so $4^\eta$ cells for the $x$-type variables over each cell of the simplex CAD. \halmos
\endproof

\proof{Proof of Corollary \ref{cor:cad_size_comparison}}
    The number of variables in our system is $\eta + \sum_1^\phi \tau_i$, the number of polynomials is $\phi + 1$ (the $\phi$ simplices plus $f^*$). We let the maximum degree (over all variables, over all polynomials) be denoted as $d$: note that in this corollary we assume that $f^*$ is linear in the $x$-type variables, but we make no assumptions on the degree of the $\alpha$-type variables, so we simply denote it as $d$. Based on the complexity analysis in \citet{england2015improving}, the dominant term of the bound of the number of cells of a general CAD, $N$, substituting our values, is:
    \begin{equation*}
        (2d)^{2^{\eta + \sum_1^\phi \tau_i} - 1} (\phi+1)^{2^{\eta + \sum_1^\phi \tau_i} - 1} 2^{2^{\eta + \sum_1^\phi \tau_i - 1} - 1}
    \end{equation*}

    Whereas when solving with the methods we have developed, we have the bound on $N_M$:
    \begin{align*}
        3^{\sum_{i=1}^\phi \tau_i} \times 4^\eta <
        4^{\eta + \sum_{i=1}^\phi \tau_i}
    \end{align*}

    For simplicity of notation, let $\psi = \eta + \sum_{i=1}^\phi \tau_i$. Then, we will prove that $4^\psi = o(2^{2^{\psi - 1} - 1})$. Consider the following ratio:
    \begin{align*}
        &\frac{4^\psi}{2^{2^{\psi - 1} - 1}} = \frac{2^{2\psi}}{2^{2^{\psi - 1} - 1}} = 2^{2\psi - 2^{\psi - 1} + 1}
    \end{align*}

    As $\psi \to \infty$, the exponent becomes large and negative, so that the term itself goes to $0$. So, $4^\psi = o(2^{2^{\psi - 1} - 1})$. Hence, by properties of little-o, $N_M = o(N)$. \halmos
\endproof

\proof{Proof of Theorem \ref{thm:nphard}}
    Suppose we have an instance of 3-SAT, i.e., Boolean satisfiability with three literals in each clause, with $m$ clauses. We refer the reader to \cite{papa} for full details on the definition of 3-SAT. We first exhibit a reduction to a simplex-extensible instance of system \eqref{eq:general_form} with $f^* = g_0 + \sum_{i=1}^\eta$. 
    
    For each Boolean variable $z_i$ the 3-SAT instance, introduce variables $\alpha_{i,1}$ (representing $z_i$ being true) and $\alpha_{i,2}$ (representing $z_i$ being false), with the constraints that $0 \leq \alpha_{i,1}, \alpha_{i,2} \leq 1$ and $\alpha_{i,1} + \alpha_{i,2} = 1$. Note that this means that in the context of system \eqref{eq:general_form}, we have $\phi$ equalling the number of Boolean variables, and $\tau = 2$ for each simplex. With these $\alpha$ variables, form the polynomial $f_1 = \sum_{i=1}^\phi \sum_{j=1}^2 \alpha_{i,j} (1 - \alpha_{i,j})$. Observe that $f_1 = 0$ if and only if all of the $\alpha$-type variables are either $0$ or $1$.

    Next, take a clause $C_i = (\ell_{i,1}, \ell_{i,2}, \ell_{i,3})$ in the 3-SAT instance. Define the term $I_{i,j}$ which is $1-\alpha_{k,1}$ if $\ell_{i,j} = z_k$ or $\alpha_{k,1}$ if $\ell_{i,j} = \neg z_k$. So, for each clause $C_i$, form the product $I_{i,1} I_{i,2} I_{i,3}$, which is clearly $0$ if at least one of the constituent $\alpha$ variables is 0, i.e., if at least one of the respective $z$ variables is 1 -- in other words, it is $0$ if and only if the clause $C_i$ is true. Then, we take the sum over all clauses, forming the polynomial $f_2 = \sum_{i=1}^m I_{i,1} I_{i,2} I_{i,3}$. Hence, $f_2 = 0$ if and only if all clauses are true.

    Finally, form the following polynomial inequality:

    \begin{align*}
        f^* := -\left( \underbrace{\sum_{i=1}^\phi \sum_{j=1}^2 \alpha_{i,j}}_{f_1} + \underbrace{\sum_{i=1}^m I_{i,1} I_{i,2} I_{i,3}}_{f_2} \right) \geq 0
    \end{align*}

    We now have an instance of system \eqref{eq:general_form}. We may observe that there are no $x$-type variables here: we can vacuously add them to $f^*$ as, e.g., by adding $x-x$; this will have no effect on our construction. All the $\alpha$ variables are continuous between $0$ and $1$, and lie on disjoint simplices. This $f^*$ is indeed in the form $g_0 + \sum_{i=1}^0 g_i x_i$. As there are no $x$ terms, we can simply set $g_0 = f^*$. Furthermore, as both $f_1$ and $f_2$ are non-negative, then $g_0$ is sign-invariant. Hence, this system is simplex-extensible. 

    We now have to show that this reduction is valid. First, suppose we have a satisfying assignment of the $z$ variables to the 3-SAT instance. This implies that all of the $\alpha$ variables are either $0$ or $1$, by construction, so $f_1=0$, and due to all clauses being true we have $f_2=0$. So, $f^* = -(0 + 0) = 0 \geq 0$. So, the instance of system \eqref{eq:general_form} is feasible.

    Secondly, suppose we have a feasible solution to the instance of system \eqref{eq:general_form}. Since $f_1, f_2 \geq 0$ always, we have that $f^* \leq 0$ always. So, with a feasible solution, we have $f^* \leq 0$ and $f^* \geq 0$, so $f^* = 0$. This implies that both $f_1 = 0$ and $f_2 = 0$. If $f_1 = 0$, all the $\alpha$ variables are either $0$ or $1$, meaning we can form valid Boolean variables $z$ from them. Secondly, if $f_2 = 0$, we also have a satisfying assignment to the 3-SAT instance.

    Hence, the constructed instance of \eqref{eq:general_form} is feasible if and only if the original instance of 3-SAT is satisfiable. The NP-hardness of simplex-extensible \eqref{eq:general_form} with $f^* = g_0 + \sum_{i=1}^\eta g_i$ and the NP-hardness of general \eqref{eq:general_form} follow. \halmos
\endproof

\begin{remark}
    Note that the construction of $f_2$ uses very common techniques that are used to show complexity of polynomial systems. In an arbitrary polynomial system, we can fix the $\alpha$ variables to be binary by introducing the constraint $\alpha(\alpha-1)=0$. However, we cannot introduce additional constraints when forming system \eqref{eq:general_form}. We get around this by incorporating the constraints implicitly into $f^*$, so that the $\alpha$ still end up being constrained to be binary.
\end{remark}

\proof{Proof of Corollary \ref{cor:simplex_ext_f_nonneg}}
    If a polynomial has no non-negative roots, it is sign-invariant over the non-negative reals. So, for each $i$, substitute $z_i = f_i$ if $f_i$ is positive or $z_i = -f_i$ if $f_i$ is negative, and treat as a system with $z_i$, noting that $z_i \geq 0$. The simplex-extensibility follows from Theorem \ref{thm:simplex_ext_xg} and the projection-definability follows from Lemma \ref{lem:xg_projdef}. Then, build the CAD and substitute $f_i$ for $z_i$ at the end. \halmos
\endproof

\section{Geometry of a two-way sensitivity analysis}
\label{sec:twowaygeom}

The formulations of the total reward inequalities allow us to formalize some facts about the geometry of a two-way sensitivity analysis in a general Markov reward process. Suppose we assert $R_\infty \geq T$, which can be reformulated as the inequality $f^* \geq 0$, where $f^*$ is a multilinear polynomial in the Markov chain's parameters. Since $f^*$ is multilinear, if we fix all but two of the parameters, what remains is a polynomial that is either linear (i.e., no variables being multiplied together) or bilinear (i.e., contains a term where both variables are being multiplied together). If it is linear, then the solution to the inequality is a half-plane. If it is bilinear, then the boundary of the solution to the inequality is a hyperbola, and the valid space is one of the sides of this hyperbola. Depending on the shape of the hyperbola and which side the solution lies -- which depends on the values of the fixed parameters -- the solution is either convex or concave. 

We can discuss these cases more specifically by analyzing which types of parameters are chosen to vary. If both are rewards, then $f^*$ is a linear function, and so the valid parameter space is a half-plane. If one variable is a reward, and another is an entry in $\mat{\pi}$ or $\mat{P}$, then $f^*$ is (possibly) bilinear, due to the two variables being multiplied together. Therefore, the boundary of the valid parameter space is a rational function, and the valid parameter space lies below or above this space. Hence, depending on the boundary and which side of the boundary is valid, we can get a convex region or a concave region. Note that if there is only a single free entry in $\mat{\pi}$ or $\mat{P}$, then it is uniquely identified by the other entries due to the stochastic constraints, so this would reduce to a one-way sensitivity analysis. Hence, for a two-sensitivity analysis to be meaningful, we want this free parameter to be bound by an inequality, e.g., if we have $\pi = [x, 1-x]$. 

Lastly, we consider the case where both free parameters are elements of $\mat{\pi}$ or $\mat{P}$. Firstly note that if both free parameters lie on the same simplex, i.e., both are from $\mat{\pi}$ or both are from the same row of $\mat{P}$, then by simple substitution we can rewrite $f^*$ as a linear function in only one of the variables, and solve for its range. Then, on the two-dimensional plane, the set of valid parameter values is in fact only a line segment. On the other hand, if they do not lie on the same simplex, then $f^*$ is bilinear or linear, and the above arguments about the geometry follow. For example, one case where $f^*$ is definitely linear is when the two free parameters are elements of $\mat{P}$ that are in the same column. Then, for the entry in the adjugate of $\mat{I} - \lambda \mat{P}$ corresponding to one of the variables, the other one will be absent. A similar argument, based on the Laplace expansion definition of the determinant, applies to $\det(\mat{I} - \lambda \mat{P})$ that the two variables will not be multiplied together. We summarize these cases in Table \ref{tab:twoway_cases}.

\begin{table}[]
    \centering
    \begin{tabular}{p{5cm}|c|c}
        \toprule
        Free parameters & Form of $f^*$ & Valid parameter space \\
        \midrule
        Both from $\mat{r}$ & Bivariate linear & Half-plane \\
        One from $\mat{r}$, one from $\mat{\pi}$ or $\mat{P}$ & Bilinear or bivariate linear & One side of hyperbola or half-plane \\
        Both from $\mat{\pi}$ or $\mat{P}$ & & \\
        --- On same simplex & Univariate linear & Line segment \\
        --- Same column of $\mat{P}$ & Bivariate linear & Half-plane \\
        --- Otherwise & Bilinear or bivariate linear & One side of hyperbola or half-plane \\
        \bottomrule
    \end{tabular}
    \caption{Geometry of two-way sensitivity analysis of the infinite horizon reward of a Markov reward process}
    \label{tab:twoway_cases}
\end{table}

The geometry of the \textit{finite horizon} parameter space is almost the same, as the function $f^*$ is still linear in the entries of $\mat{r}$ and linear in the entries of $\mat{\pi}$. However, due to summing exponents of $\mat{P}$, the entries of $\mat{P}$ appear with exponents (i.e., with degrees higher than 1, unless the reward is only computed for a single period). So, if an entry from $\mat{P}$ is chosen, we possibly obtain a polynomial with degree greater than one -- not necessarily a bilinear function as in the infinite horizon case. In this case, it is difficult to solve exactly for one variable in terms of another, as to do so may require an arduous expression with radicals. And indeed, it may be impossible if we sum the rewards of the Markov process for 5 periods or more, due to the Abel-Ruffini Theorem \citep{ruffini, abel}. However, the boundary is still always either convex or concave, and depending on the shape of the curve and which side of the curve we are on, the valid parameter space may be non-convex.

\end{document}